\documentclass[final]{elsart}

\usepackage{longtable}
\usepackage{graphicx}
\usepackage{natbib}

\begin{document}

\bibliographystyle{elsart-harv}

\begin{frontmatter}

\title{Square root voting in the Council of the European Union: Rounding effects and the Jagiellonian Compromise}

\author{Martin Kurth}

\ead{martin.kurth@nottingham.ac.uk}

\address{School of Mathematical Sciences, University of Nottingham, University Park, Nottingham, NG7 2RD, UK}

\begin{abstract}
In recent years, enlargement of the European Union has brought with it renewed discussion
of voting arrangements in the Council of the EU. During these negotiations, the Polish
government proposed a voting scheme that gives each country a voting weight proportional
to the square root of its population, and sets a quota according to an optimality condition
(``Jagiellonian Compromise''). In this paper, the optimal quota is found exactly for the
current population data from the 27 EU member states, and it is found that rounding of the
voting weights can be used to improve the voting scheme.
\end{abstract}

\end{frontmatter}

\section{Introduction}

Following the failure of the draft EU constitution in referenda in France and the
Netherlands, and the subsequent negotiations on a Reform Treaty, the voting arrangements
in the Council of the EU have received a significant amount of public attention, not least
through the Polish proposal of a voting system that gives every member state a voting
weight proportional to the square root of its population, \citet{euobserver}. The idea of
this voting scheme is to give every EU citizen the same influence on decisions in the
Council, based on an analysis of their voting power. The concept of voting power used
in this context was first introduced by \citet{penrose1}, and adapted to
the EU framework by \citet{slomczynski1}, \citet{slomczynski3}.

The purpose of this paper is twofold: To find an exact value for the optimal quota
required to make the voting scheme as close as possible to the ideal of equal influence
for all citizens, and to demonstrate the impact that rounding effects have on this
quota.

The paper is organised as follows: In section~\ref{s_squarerootvoting}, the main
ideas behind square root voting are summarised, and some notation is
introduced. Section~\ref{s_rounding} examines the effect of using rounded rather
than exact voting weights by analysing two example scenarios. In the example in
subsection~\ref{ss_example1}, the rounding makes the voting system less fair, while
the example in subsection~\ref{ss_example2} demonstrates that rounding can actually
improve the voting scheme. Section~\ref{s_euvoting} is devoted to voting in the
Council of the EU. The current voting arrangements are briefly described in
subsection~\ref{ss_current}, while square root voting is applied to the EU
context in subsection~\ref{ss_eusquareroot} following \citet{slomczynski1}, and rounding effects are analysed. A similar
analysis including member quotas is carried out in subsection~\ref{ss_eumemberquota}.
Section~\ref{s_summary} summarises the results.

\section{Voting power and square root voting}
\label{s_squarerootvoting}

The concept of voting power was first introduced by \citet{penrose1} for
the UN General Assembly, but his analysis applies to all bodies where votes are cast
by members elected in separate geographical areas with different numbers of voters,
and the elected members from each area cast a joint vote. This is, for example, the case in
the Electoral College in US Presidential elections, and in the Council of the EU,
where the national governments of the (currently) 27 member states are represented.

As the numbers of citizens in the 27 member states differ hugely (between approximately
400,000 in Malta and 82.3 million in Germany), the task at hand is to create a ``fair''
voting system by assigning a voting weight to each country in the Council, where the
most reasonable definition of ``fair'' is that every voter should have the same influence
on the outcome of votes in the Council. Thus a mathematically rigorous definition of
``influence'' is needed.

The first step is to analyse the voting power of each voter within their own country.
Penrose defines voting power as the probability for the vote of one individual voter
to change the voting outcome, assuming that all other voters vote randomly.

In the following discussion, let us assume that there are only two voting options,
``yes'' and ``no'', and that the outcome of the vote is positive if there are more
``yes'' votes than ``no'' votes, and negative otherwise.

In a country with $N$ voters, the total number of possible results of the vote is
$2^N$. ``Random voting'' means that all these results are equally likely.

If $N$ is odd, the probability for a single voter to make a difference is the probability
that the remaining $N-1$ votes are evenly split between ``yes'' and ``no'' votes,
which is
\begin{equation}
  P_{\rm odd}=\frac{1}{2^{N-1}}\left(\begin{array}{c}N-1\\(N-1)/2\end{array}\right).
\end{equation}
If $N$ is even, the probability for a single vote to be decisive is the probability that
without this vote, there is one ``yes'' vote more than there are ``no'' votes, which
means
\begin{equation}
  P_{\rm even}=\frac{1}{2^{N-1}}\left(\begin{array}{c}N-1\\N/2\end{array}\right).
\end{equation}
Using the Stirling approximation for factorials, and the fact that
\begin{equation}
  \sqrt{N-1}\approx\sqrt{N} \qquad (N\gg 1),
\end{equation}
we obtain for the probability of a single vote to be decisive (ie the voting power
of an individual voter)
\begin{equation}
  P_N \propto \frac{1}{\sqrt{N}}.
\end{equation}
To give every voter the same influence on the voting outcome in the Council of the
EU, voting weights there must be assigned in such a way, that the voting power
of each country is proportional to the square root of the number of its inhabitants.
This is known as \emph{Penrose's Square Root Law}.

The analysis of voting power within the Council (or other body) is also due to Penrose,
but was independently carried out by \cite{banzhaf1}.
Let us assume the number of members in the decision taking body is $n$ (in the case
of the Council of the EU $n=27$). All countries will be assigned voting weights, and
in the following the voting weight of country $i$ will be denoted by $v_i$, where the
weights are normalised such that
\begin{equation}
\label{weight_normalisation}
  \sum_{i=1}^n v_i=1.
\end{equation}
Any subset of
countries is called a coalition $C$, and there will be a quota $R$ which has to be
reached for a positive voting outcome, ie a coalition $C$ will be successful if
\begin{equation}
  \sum_{i\in C}v_i\geq R.
\end{equation}
The voting power of country $i$ is defined as the probability that its vote is
critical (ie makes a difference to the voting outcome) under the assumption that
the remaining $n-1$ countries vote randomly. If the total number of winning coalitions
(ie coalitions that reach the quota) is $\omega$, and the number of winning coalitions
that contain country $i$ is $\omega_i$, there are
\begin{equation}
  \label{eta_definition}
  \eta_i=2\omega_i-\omega
\end{equation}
coalitions where the vote of country $i$ is critical, ie which would not reach
the quota without country $i$'s vote. As the vote of the remaining $n-1$ countries
can yield $2^{n-1}$ different outcomes, which are all equally likely, the probability
for country $i$'s vote to be critical is
\begin{equation}
\label{absolute_banzhaf_definition}
  B_i=\frac{\eta_i}{2^{n-1}}.
\end{equation}
The number $B_i$ is called the absolute Banzhaf index. For our purposes, it is useful
to introduce the relative Banzhaf index
\begin{equation}
\label{relative_banzhaf_definition}
  \beta_i=\frac{B_i}{\sum_{j=1}^n B_j}=\frac{\eta_i}{\sum_{j=1}^n \eta_j}.
\end{equation}
To make the voting system fair, we would ideally want
\begin{equation}
  \beta_i=\frac{\sqrt{N_i}}{\sum_{j=1}^n \sqrt{N_j}}
\end{equation}
for all countries $i$, where $N_i$ is the population of country $i$. Such weights may not
exist, and even finding an approximation by minimising an appropriately defined deviation of the
voting powers from their ideal values is a complicated multi-dimensional optimisation
problem. To approximate the ideal voting power distribution, Penrose suggested to use
\emph{voting weights} proportional to the square root of the respective populations, ie
\begin{equation}
\label{weight_definition}
  v_i=\frac{\sqrt{N_i}}{\sum_{j=1}^n \sqrt{N_j}}
\end{equation}
for all countries $i$. This
is the scheme that \cite{slomczynski1} and \cite{slomczynski3} propose for the Council of the EU, where they choose
the quota $R$ such that
\begin{equation}
\label{sigma_definition}
  \sigma=\sqrt{\frac{1}{n}\sum_{i=1}^n\left(v_i-\beta_i\right)^2}
\end{equation}
is minimal. The voting scheme constructed in this way is called the Jagiellonian Compromise.

A further useful notion is the \emph{efficiency} $\epsilon$ of a voting system (also known as the
\emph{Coleman index}), which is the probability
for a randomly chosen coalition to be successful,
\begin{equation}
  \epsilon = \frac{\omega}{2^n}.
\end{equation}

The concept of square root voting power has been criticised by many authors (see \cite{Garrett1}, \cite{Albert1}, \cite{bafumi1} for
recent criticism)
on the grounds that actual voter behaviour is far from random. In \cite{bafumi1}, it is shown that
in US presidential elections, the random voting assumption seems to overestimate the probability of close voting outcomes in
larger states, and
the authors contrast \emph{actual voting power} with the random voting concept. Slomczy{\'n}ski
and Zyczkowski, as well as for example \cite{Leech2} (following \cite{Shapley1}, \cite{banzhaf2}, \cite{Coleman1}) point out that
the voting power concept derived from the random voting assumption is an \emph{a priori voting power}, based
on the freedom of the voters to make their choice in any way they want, rather than on the way they actually
vote. It is hard to see how a voting system could be based on actual voting power, as this would have to be
determined from previous election results. This would mean that voters could influence the weight of their
votes in future elections by the way they vote. A further problem with the idea of basing voting weights on results of earlier
elections is that voter behaviour can change over time, and older election results may
not reflect the actual voting power held by current voters. Actual voting power in Europe may
be quite different from voting power in America. Since January 1997, the governments of all 27 EU member states have changed
at least once as a result of a general election, which may indicate that close election outcomes are more
common in European countries than in the larger US states. Determining actual voting power would be a very ambitious
project, involving very large separate data sets from all EU member states, including a large number of election
results and/or sufficiently accurate opinion polls to derive statistical models for voter behaviour in all
countries currently in the EU. In the absence of such data, and for the reason of making voting weights independent
from how voters have voted in the past, the square root voting system seems to be the best one can do. At least
it is based on clearly stated assumptions, that can in principle be tested empirically.
This is not the case for the current EU voting arrangements.

The discussion on what definition of voting power to use is ongoing, and is not the subject of this paper. Here we will
regard a voting system based on the Penrose definition of voting power as ``fair''.

It can reasonably be argued that voting within the Council is not random either: Governments will, for example, form
strategic alliances. However, such voting behaviour is hard to predict and to include in a model of actual voting
power, and will not be considered here (see~\cite{Hosli1} for a more detailed discussion).

\section{Rounding errors}
\label{s_rounding}

The voting weights defined in equation~\ref{weight_definition} will in general not be rational numbers, and some
rounding rule will have to be applied; even if no precise rounding rule is given, the weights will be implicitly
rounded by the accuracy of the floating point representation on the computer system used. In \cite{slomczynski1},
the authors
round all voting weights to four decimal places behind the decimal point (although it is not entirely clear if
their minimisation of $\sigma$ is carried out with these rounded weights, or if some more accurate values are
used in the calculation). The purpose of this section is to show that the rounding can have a very large effect
on the optimal quota $R$ and on the minimal value of $\sigma$ that can be achieved by the Penrose square root
voting system. In the following, ``$k$ digit rounding'' means standard rounding of voting weights to $k$ decimal
places, where voting weights are normalised according to equation~\ref{weight_normalisation}.

\subsection{Example 1}
\label{ss_example1}

To demonstrate the effect of rounding errors on the minimisation of $\sigma$, let us consider an example with
4 countries, and different numbers of voters, as shown in table~\ref{t_ex1_countries}.

\begin{table}[ht]
\begin{center}
\begin{tabular}{|c|r|r|r|r|}
\hline
Country & \multicolumn{1}{c|}{$N_i$} & \multicolumn{1}{c|}{$\sqrt{N_i}$} & \multicolumn{1}{c|}{$v_i$}& \multicolumn{1}{c|}{$v_i$ (3 digits)}\\
\hline
1 & 20,295,025 & 4505 & 0.4505 & 0.451 \\
2 & 6,265,009 & 2503 & 0.2503 & 0.250 \\
3 & 4,012,009 & 2003 & 0.2003 & 0.200 \\
4 & 978,121 & 989 & 0.0989 & 0.099 \\
\hline
\end{tabular}
\end{center}
\caption{Example 1: Populations and voting weights for 4 countries.}
\label{t_ex1_countries}
\end{table}

In this example, 4 digit rounding is exact. Using the 4 digit voting weights, we can list the voting weights for all
possible coalitions, see table~\ref{t_ex1_coalitions4}.

\begin{table}[ht]
\begin{center}
\begin{tabular}{|c|r||c|r|}
\hline
Countries & $\sum v_i$ & Countries & $\sum v_i$ \\
\hline
1, 2, 3, 4 & 1.0000 & 2, 3 & 0.4506 \\
1, 2, 3 & 0.9011 & 1 & 0.4505 \\
1, 2, 4 & 0.7997 & 2, 4 & 0.3492 \\
1, 3, 4 & 0.7497 & 3, 4 & 0.2992 \\
1, 2 & 0.7008 & 2 & 0.2503 \\
1, 3 & 0.6508 & 3 & 0.2003 \\
2, 3, 4 & 0.5495 & 4 & 0.0989 \\
1, 4 & 0.5494 & None  & 0.0000 \\
\hline
\end{tabular}
\end{center}
\caption{Example 1: Coalitions and their voting weights without rounding.}
\label{t_ex1_coalitions4}
\end{table}

From this table, we can read off the number $\omega$ of coalitions having a majority for any given quota $R$,
as well as the numbers $\omega_i$ of these coalitions containing country $i$. This means that for all quotas
$0\leq R\leq 1$ and all countries $i$, the numbers $\eta_i$ can be calculated from equation~\ref{eta_definition},
giving the Banzhaf index, and ultimately the deviation $\sigma$ from a fair square root voting system. For example,
for $0.9011<R\leq 1$, there is exactly one successful coalition, containing all four countries, which
means that $\eta_i=1$ for $i=1,\ldots,4$. The deviation $\sigma$ then follows from equations~\ref{absolute_banzhaf_definition},
\ref{relative_banzhaf_definition} and \ref{sigma_definition}. The values of $\sigma$ for all quotas $R$ are shown
in table~\ref{t_ex1_deviation4}.

\begin{table}[ht]
\begin{center}
\begin{tabular}{|r|r||c|c|c|c||r|}
\hline
\multicolumn{1}{|c|}{$R>$} & \multicolumn{1}{c||}{$R\leq$} & $\eta_1$ & $\eta_2$
 & $\eta_3$ & $\eta_4$ & \multicolumn{1}{c|}{$\sigma$} \\
\hline
 $0.9011$ & $1.0000$ & $1$ & $1$ & $1$ & $1$ & $0.127966$ \\
 $0.7997$ & $0.9011$ & $2$ & $2$ & $2$ & $0$ & $0.109660$ \\
 $0.7497$ & $0.7997$ & $3$ & $3$ & $1$ & $1$ & $0.083069$ \\
 $0.7008$ & $0.7497$ & $4$ & $2$ & $2$ & $2$ & $0.061850$ \\
 $0.6508$ & $0.7008$ & $5$ & $3$ & $1$ & $1$ & $0.061200$ \\
 $0.5495$ & $0.6508$ & $6$ & $2$ & $2$ & $0$ & $0.093088$ \\
 $0.5494$ & $0.5495$ & $5$ & $3$ & $3$ & $1$ & $0.031053$ \\
 $0.4506$ & $0.5494$ & $6$ & $2$ & $2$ & $2$ & $0.061580$ \\
 $0.4505$ & $0.4506$ & $5$ & $3$ & $3$ & $1$ & $0.031053$ \\
 $0.3492$ & $0.4505$ & $6$ & $2$ & $2$ & $0$ & $0.093088$ \\
 $0.2992$ & $0.3492$ & $5$ & $3$ & $1$ & $1$ & $0.061200$ \\
 $0.2503$ & $0.2992$ & $4$ & $2$ & $2$ & $2$ & $0.061850$ \\
 $0.2003$ & $0.2503$ & $3$ & $3$ & $1$ & $1$ & $0.083069$ \\
 $0.0989$ & $0.2003$ & $2$ & $2$ & $2$ & $0$ & $0.109660$ \\
 $0.0000$ & $0.0989$ & $1$ & $1$ & $1$ & $1$ & $0.127966$ \\
\hline
\end{tabular}
\end{center}
\caption{Example 1: Deviation from fair voting without rounding.}
\label{t_ex1_deviation4}
\end{table}

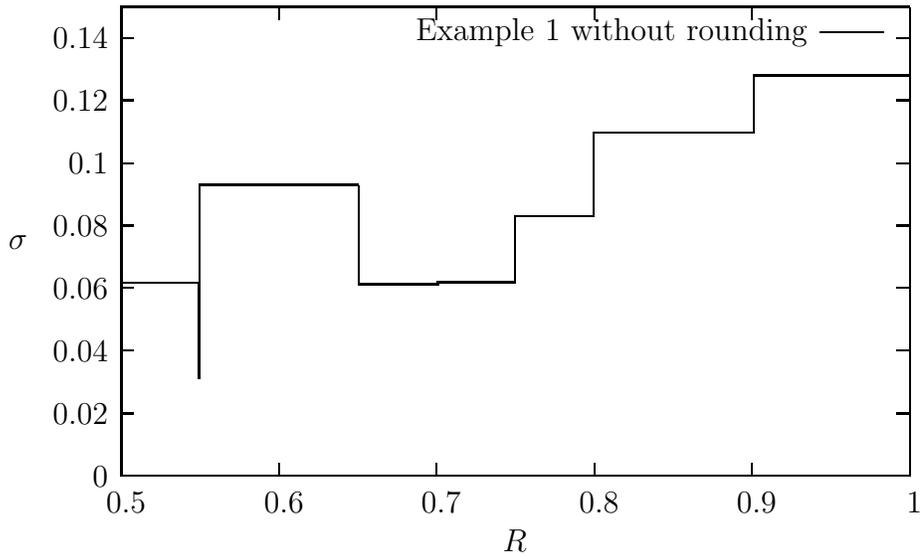
\begin{figure}[ht]
\begin{center}
% GNUPLOT: LaTeX picture
\setlength{\unitlength}{0.240900pt}
\ifx\plotpoint\undefined\newsavebox{\plotpoint}\fi
\sbox{\plotpoint}{\rule[-0.200pt]{0.400pt}{0.400pt}}%
\begin{picture}(1500,900)(0,0)
\sbox{\plotpoint}{\rule[-0.200pt]{0.400pt}{0.400pt}}%
\put(201.0,123.0){\rule[-0.200pt]{4.818pt}{0.400pt}}
\put(181,123){\makebox(0,0)[r]{ 0}}
\put(1419.0,123.0){\rule[-0.200pt]{4.818pt}{0.400pt}}
\put(201.0,221.0){\rule[-0.200pt]{4.818pt}{0.400pt}}
\put(181,221){\makebox(0,0)[r]{ 0.02}}
\put(1419.0,221.0){\rule[-0.200pt]{4.818pt}{0.400pt}}
\put(201.0,320.0){\rule[-0.200pt]{4.818pt}{0.400pt}}
\put(181,320){\makebox(0,0)[r]{ 0.04}}
\put(1419.0,320.0){\rule[-0.200pt]{4.818pt}{0.400pt}}
\put(201.0,418.0){\rule[-0.200pt]{4.818pt}{0.400pt}}
\put(181,418){\makebox(0,0)[r]{ 0.06}}
\put(1419.0,418.0){\rule[-0.200pt]{4.818pt}{0.400pt}}
\put(201.0,516.0){\rule[-0.200pt]{4.818pt}{0.400pt}}
\put(181,516){\makebox(0,0)[r]{ 0.08}}
\put(1419.0,516.0){\rule[-0.200pt]{4.818pt}{0.400pt}}
\put(201.0,614.0){\rule[-0.200pt]{4.818pt}{0.400pt}}
\put(181,614){\makebox(0,0)[r]{ 0.1}}
\put(1419.0,614.0){\rule[-0.200pt]{4.818pt}{0.400pt}}
\put(201.0,713.0){\rule[-0.200pt]{4.818pt}{0.400pt}}
\put(181,713){\makebox(0,0)[r]{ 0.12}}
\put(1419.0,713.0){\rule[-0.200pt]{4.818pt}{0.400pt}}
\put(201.0,811.0){\rule[-0.200pt]{4.818pt}{0.400pt}}
\put(181,811){\makebox(0,0)[r]{ 0.14}}
\put(1419.0,811.0){\rule[-0.200pt]{4.818pt}{0.400pt}}
\put(201.0,123.0){\rule[-0.200pt]{0.400pt}{4.818pt}}
\put(201,82){\makebox(0,0){ 0.5}}
\put(201.0,840.0){\rule[-0.200pt]{0.400pt}{4.818pt}}
\put(449.0,123.0){\rule[-0.200pt]{0.400pt}{4.818pt}}
\put(449,82){\makebox(0,0){ 0.6}}
\put(449.0,840.0){\rule[-0.200pt]{0.400pt}{4.818pt}}
\put(696.0,123.0){\rule[-0.200pt]{0.400pt}{4.818pt}}
\put(696,82){\makebox(0,0){ 0.7}}
\put(696.0,840.0){\rule[-0.200pt]{0.400pt}{4.818pt}}
\put(944.0,123.0){\rule[-0.200pt]{0.400pt}{4.818pt}}
\put(944,82){\makebox(0,0){ 0.8}}
\put(944.0,840.0){\rule[-0.200pt]{0.400pt}{4.818pt}}
\put(1191.0,123.0){\rule[-0.200pt]{0.400pt}{4.818pt}}
\put(1191,82){\makebox(0,0){ 0.9}}
\put(1191.0,840.0){\rule[-0.200pt]{0.400pt}{4.818pt}}
\put(1439.0,123.0){\rule[-0.200pt]{0.400pt}{4.818pt}}
\put(1439,82){\makebox(0,0){ 1}}
\put(1439.0,840.0){\rule[-0.200pt]{0.400pt}{4.818pt}}
\put(201.0,123.0){\rule[-0.200pt]{298.234pt}{0.400pt}}
\put(1439.0,123.0){\rule[-0.200pt]{0.400pt}{177.543pt}}
\put(201.0,860.0){\rule[-0.200pt]{298.234pt}{0.400pt}}
\put(201.0,123.0){\rule[-0.200pt]{0.400pt}{177.543pt}}
\put(40,491){\makebox(0,0){$\sigma$}}
\put(820,21){\makebox(0,0){$R$}}
\put(1279,820){\makebox(0,0)[r]{Example 1 without rounding}}
\put(1299.0,820.0){\rule[-0.200pt]{24.090pt}{0.400pt}}
\put(201,426){\usebox{\plotpoint}}
\put(201.0,426.0){\rule[-0.200pt]{29.390pt}{0.400pt}}
\put(323.0,276.0){\rule[-0.200pt]{0.400pt}{36.135pt}}
\put(323.0,276.0){\usebox{\plotpoint}}
\put(324.0,276.0){\rule[-0.200pt]{0.400pt}{73.234pt}}
\put(324.0,580.0){\rule[-0.200pt]{60.225pt}{0.400pt}}
\put(574.0,424.0){\rule[-0.200pt]{0.400pt}{37.580pt}}
\put(574.0,424.0){\rule[-0.200pt]{29.872pt}{0.400pt}}
\put(698.0,424.0){\rule[-0.200pt]{0.400pt}{0.723pt}}
\put(698.0,427.0){\rule[-0.200pt]{29.149pt}{0.400pt}}
\put(819.0,427.0){\rule[-0.200pt]{0.400pt}{25.054pt}}
\put(819.0,531.0){\rule[-0.200pt]{29.872pt}{0.400pt}}
\put(943.0,531.0){\rule[-0.200pt]{0.400pt}{31.558pt}}
\put(943.0,662.0){\rule[-0.200pt]{60.466pt}{0.400pt}}
\put(1194.0,662.0){\rule[-0.200pt]{0.400pt}{21.681pt}}
\put(1194.0,752.0){\rule[-0.200pt]{59.020pt}{0.400pt}}
\put(201.0,123.0){\rule[-0.200pt]{298.234pt}{0.400pt}}
\put(1439.0,123.0){\rule[-0.200pt]{0.400pt}{177.543pt}}
\put(201.0,860.0){\rule[-0.200pt]{298.234pt}{0.400pt}}
\put(201.0,123.0){\rule[-0.200pt]{0.400pt}{177.543pt}}
\end{picture}
\end{center}
\caption{Example 1: Deviation from fair voting without rounding.}
\label{f_ex1_deviation4}
\end{figure}

A plot of $\sigma$ as a function of the quota $R$ is shown in figure~\ref{f_ex1_deviation4}. Only
quotas $R>=0.5$ are shown here, as quotas less than $0.5$ would mean that a coalition could be
successful even if the opposing coalition has more votes.

The minimum of $\sigma$ can be read off. We have
\begin{equation}
\sigma_{\rm min}=0.031053,
\end{equation}
and the minimum occurs in the interval
\begin{equation}
0.5494<R\leq 0.5495.
\end{equation}

The situation changes dramatically if $\sigma$ is calculated from the rounded voting weights. For this case,
the possible coalitions and their voting weights are shown in table~\ref{t_ex1_coalitions3}.

\begin{table}[ht]
\begin{center}
\begin{tabular}{|c|r||c|r|}
\hline
Countries & $\sum v_i$ & Countries & $\sum v_i$ \\
\hline
1, 2, 3, 4 & 1.000 & 1 & 0.451 \\
1, 2, 3 & 0.901 & 2, 3 & 0.450 \\
1, 2, 4 & 0.800 & 2, 4 & 0.349 \\
1, 3, 4 & 0.750 & 3, 4 & 0.299 \\
1, 2 & 0.701 & 2 & 0.250 \\
1, 3 & 0.651 & 3 & 0.200 \\
1, 4 & 0.550 & 4 & 0.099 \\
2, 3, 4 & 0.549 & None & 0.000 \\
\hline
\end{tabular}
\end{center}
\caption{Example 1: Coalitions and their voting weights with 3 digit rounding.}
\label{t_ex1_coalitions3}
\end{table}

Comparison with table~\ref{t_ex1_coalitions4} shows that the coalition of countries 2, 3 and 4, and the coalition
of countries 1 and 4 have changed places. This leads to different results for $\eta_i$ and $\sigma$, as shown in
table~\ref{t_ex1_deviation3}.

\begin{table}[ht]
\begin{center}
\begin{tabular}{|r|r||c|c|c|c||r|}
\hline
\multicolumn{1}{|c|}{$R>$} & \multicolumn{1}{c||}{$R\leq$} & $\eta_1$ & $\eta_2$ & $\eta_3$ & $\eta_4$ & \multicolumn{1}{c|}{$\sigma$} \\ 
\hline
 $0.901$ & $1.000$ & $1$ & $1$ & $1$ & $1$ & $0.127966$ \\ 
 $0.800$ & $0.901$ & $2$ & $2$ & $2$ & $0$ & $0.109660$ \\ 
 $0.750$ & $0.800$ & $3$ & $3$ & $1$ & $1$ & $0.083069$ \\ 
 $0.701$ & $0.750$ & $4$ & $2$ & $2$ & $2$ & $0.061850$ \\ 
 $0.651$ & $0.701$ & $5$ & $3$ & $1$ & $1$ & $0.061200$ \\ 
 $0.550$ & $0.651$ & $6$ & $2$ & $2$ & $0$ & $0.093088$ \\ 
 $0.549$ & $0.550$ & $7$ & $1$ & $1$ & $1$ & $0.154031$ \\ 
 $0.451$ & $0.549$ & $6$ & $2$ & $2$ & $2$ & $0.061580$ \\ 
 $0.450$ & $0.451$ & $7$ & $1$ & $1$ & $1$ & $0.154031$ \\ 
 $0.349$ & $0.450$ & $6$ & $2$ & $2$ & $0$ & $0.093088$ \\ 
 $0.299$ & $0.349$ & $5$ & $3$ & $1$ & $1$ & $0.061200$ \\ 
 $0.250$ & $0.299$ & $4$ & $2$ & $2$ & $2$ & $0.061850$ \\ 
 $0.200$ & $0.250$ & $3$ & $3$ & $1$ & $1$ & $0.083069$ \\ 
 $0.099$ & $0.200$ & $2$ & $2$ & $2$ & $0$ & $0.109660$ \\ 
 $0.000$ & $0.099$ & $1$ & $1$ & $1$ & $1$ & $0.127966$ \\ 
\hline
\end{tabular}
\end{center}
\caption{Example 1: Deviation from fair voting with 3 digit rounding.}
\label{t_ex1_deviation3}
\end{table}

\begin{figure}[ht]
\begin{center}
% GNUPLOT: LaTeX picture
\setlength{\unitlength}{0.240900pt}
\ifx\plotpoint\undefined\newsavebox{\plotpoint}\fi
\sbox{\plotpoint}{\rule[-0.200pt]{0.400pt}{0.400pt}}%
\begin{picture}(1500,900)(0,0)
\sbox{\plotpoint}{\rule[-0.200pt]{0.400pt}{0.400pt}}%
\put(201.0,123.0){\rule[-0.200pt]{4.818pt}{0.400pt}}
\put(181,123){\makebox(0,0)[r]{ 0}}
\put(1419.0,123.0){\rule[-0.200pt]{4.818pt}{0.400pt}}
\put(201.0,221.0){\rule[-0.200pt]{4.818pt}{0.400pt}}
\put(181,221){\makebox(0,0)[r]{ 0.02}}
\put(1419.0,221.0){\rule[-0.200pt]{4.818pt}{0.400pt}}
\put(201.0,320.0){\rule[-0.200pt]{4.818pt}{0.400pt}}
\put(181,320){\makebox(0,0)[r]{ 0.04}}
\put(1419.0,320.0){\rule[-0.200pt]{4.818pt}{0.400pt}}
\put(201.0,418.0){\rule[-0.200pt]{4.818pt}{0.400pt}}
\put(181,418){\makebox(0,0)[r]{ 0.06}}
\put(1419.0,418.0){\rule[-0.200pt]{4.818pt}{0.400pt}}
\put(201.0,516.0){\rule[-0.200pt]{4.818pt}{0.400pt}}
\put(181,516){\makebox(0,0)[r]{ 0.08}}
\put(1419.0,516.0){\rule[-0.200pt]{4.818pt}{0.400pt}}
\put(201.0,614.0){\rule[-0.200pt]{4.818pt}{0.400pt}}
\put(181,614){\makebox(0,0)[r]{ 0.1}}
\put(1419.0,614.0){\rule[-0.200pt]{4.818pt}{0.400pt}}
\put(201.0,713.0){\rule[-0.200pt]{4.818pt}{0.400pt}}
\put(181,713){\makebox(0,0)[r]{ 0.12}}
\put(1419.0,713.0){\rule[-0.200pt]{4.818pt}{0.400pt}}
\put(201.0,811.0){\rule[-0.200pt]{4.818pt}{0.400pt}}
\put(181,811){\makebox(0,0)[r]{ 0.14}}
\put(1419.0,811.0){\rule[-0.200pt]{4.818pt}{0.400pt}}
\put(201.0,123.0){\rule[-0.200pt]{0.400pt}{4.818pt}}
\put(201,82){\makebox(0,0){ 0.5}}
\put(201.0,840.0){\rule[-0.200pt]{0.400pt}{4.818pt}}
\put(449.0,123.0){\rule[-0.200pt]{0.400pt}{4.818pt}}
\put(449,82){\makebox(0,0){ 0.6}}
\put(449.0,840.0){\rule[-0.200pt]{0.400pt}{4.818pt}}
\put(696.0,123.0){\rule[-0.200pt]{0.400pt}{4.818pt}}
\put(696,82){\makebox(0,0){ 0.7}}
\put(696.0,840.0){\rule[-0.200pt]{0.400pt}{4.818pt}}
\put(944.0,123.0){\rule[-0.200pt]{0.400pt}{4.818pt}}
\put(944,82){\makebox(0,0){ 0.8}}
\put(944.0,840.0){\rule[-0.200pt]{0.400pt}{4.818pt}}
\put(1191.0,123.0){\rule[-0.200pt]{0.400pt}{4.818pt}}
\put(1191,82){\makebox(0,0){ 0.9}}
\put(1191.0,840.0){\rule[-0.200pt]{0.400pt}{4.818pt}}
\put(1439.0,123.0){\rule[-0.200pt]{0.400pt}{4.818pt}}
\put(1439,82){\makebox(0,0){ 1}}
\put(1439.0,840.0){\rule[-0.200pt]{0.400pt}{4.818pt}}
\put(201.0,123.0){\rule[-0.200pt]{298.234pt}{0.400pt}}
\put(1439.0,123.0){\rule[-0.200pt]{0.400pt}{177.543pt}}
\put(201.0,860.0){\rule[-0.200pt]{298.234pt}{0.400pt}}
\put(201.0,123.0){\rule[-0.200pt]{0.400pt}{177.543pt}}
\put(40,491){\makebox(0,0){$\sigma$}}
\put(820,21){\makebox(0,0){$R$}}
\put(1279,820){\makebox(0,0)[r]{Example 1 with 3-digit rounding}}
\put(1299.0,820.0){\rule[-0.200pt]{24.090pt}{0.400pt}}
\put(201,426){\usebox{\plotpoint}}
\put(201.0,426.0){\rule[-0.200pt]{29.149pt}{0.400pt}}
\put(322.0,426.0){\rule[-0.200pt]{0.400pt}{104.551pt}}
\put(325,860){\usebox{\plotpoint}}
\put(325.0,580.0){\rule[-0.200pt]{0.400pt}{67.452pt}}
\put(325.0,580.0){\rule[-0.200pt]{60.225pt}{0.400pt}}
\put(575.0,424.0){\rule[-0.200pt]{0.400pt}{37.580pt}}
\put(575.0,424.0){\rule[-0.200pt]{29.872pt}{0.400pt}}
\put(699.0,424.0){\rule[-0.200pt]{0.400pt}{0.723pt}}
\put(699.0,427.0){\rule[-0.200pt]{29.149pt}{0.400pt}}
\put(820.0,427.0){\rule[-0.200pt]{0.400pt}{25.054pt}}
\put(820.0,531.0){\rule[-0.200pt]{29.872pt}{0.400pt}}
\put(944.0,531.0){\rule[-0.200pt]{0.400pt}{31.558pt}}
\put(944.0,662.0){\rule[-0.200pt]{60.225pt}{0.400pt}}
\put(1194.0,662.0){\rule[-0.200pt]{0.400pt}{21.681pt}}
\put(1194.0,752.0){\rule[-0.200pt]{59.020pt}{0.400pt}}
\put(201.0,123.0){\rule[-0.200pt]{298.234pt}{0.400pt}}
\put(1439.0,123.0){\rule[-0.200pt]{0.400pt}{177.543pt}}
\put(201.0,860.0){\rule[-0.200pt]{298.234pt}{0.400pt}}
\put(201.0,123.0){\rule[-0.200pt]{0.400pt}{177.543pt}}
\end{picture}
\end{center}
\caption{Example 1: Deviation from fair voting with 3 digit rounding.}
\label{f_ex1_deviation3}
\end{figure}
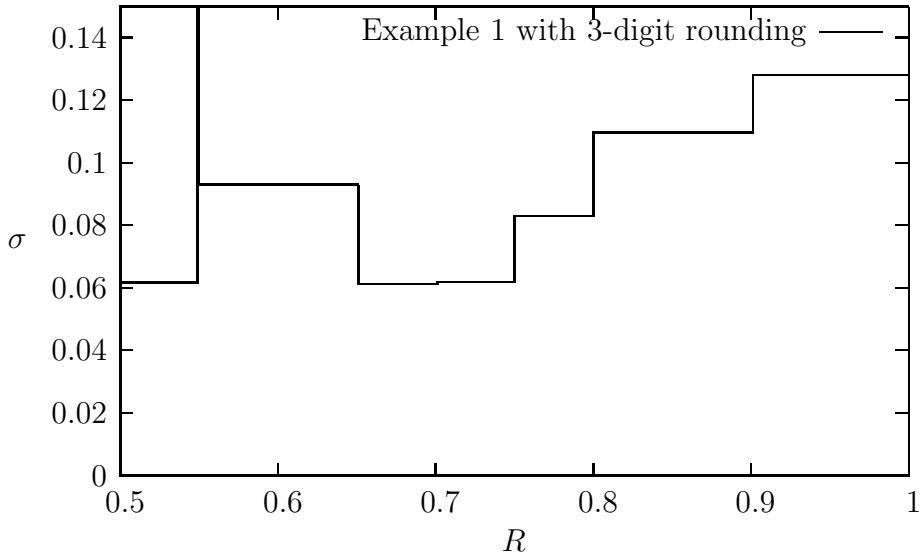

Again, $\sigma$ can be plotted as a function of $R$, as shown in figure~\ref{f_ex1_deviation3}.

Reading off the minimum of $\sigma$, we find
\begin{equation}
  \sigma_{\rm min}=0.061200,
\end{equation}
occurring at
\begin{equation}
  0.651<R\leq 0.701,
\end{equation}
which means that the minimum value of $\sigma$ has changed dramatically from the case without rounding
(it has more than doubled), and also the interval where the minimum occurs has completely changed. This shows
that rounding can have a significant impact on the minimum value of $\sigma$ and the associated quota. In this
example, rounding increases the value of $\sigma_{\rm min}$, but in general this is not necessarily the case.
By coincidence, rounding can decrease the value of $\sigma_{\rm min}$, which means that
it can improve the Jagiellonian approximation
for a voting system with square root voting power. This is illustrated in Example 2.

\subsection{Example 2}
\label{ss_example2}

This is again an example for four countries, with populations given by table~\ref{t_ex2_countries}.

\begin{table}[ht]
\begin{center}
\begin{tabular}{|c|r|r|r|r|}
\hline
Country & \multicolumn{1}{c|}{$N_i$} & \multicolumn{1}{c|}{$\sqrt{N_i}$} & \multicolumn{1}{c|}{$v_i$}& \multicolumn{1}{c|}{$v_i$ (3 digits)}\\
\hline
1 & 20,376,196 & 4514 & 0.4514 & 0.451 \\
2 & 6,280,036 & 2506 & 0.2506 & 0.251 \\
3 & 4,024,036 & 2006 & 0.2006 & 0.201 \\
4 & 948,676 & 974 & 0.0974 & 0.097 \\
\hline
\end{tabular}
\end{center}
\caption{Example 2: Populations and voting weights for 4 countries.}
\label{t_ex2_countries}
\end{table}

\begin{figure}[ht]
\begin{center}
% GNUPLOT: LaTeX picture
\setlength{\unitlength}{0.240900pt}
\ifx\plotpoint\undefined\newsavebox{\plotpoint}\fi
\sbox{\plotpoint}{\rule[-0.200pt]{0.400pt}{0.400pt}}%
\begin{picture}(1500,900)(0,0)
\sbox{\plotpoint}{\rule[-0.200pt]{0.400pt}{0.400pt}}%
\put(201.0,123.0){\rule[-0.200pt]{4.818pt}{0.400pt}}
\put(181,123){\makebox(0,0)[r]{ 0}}
\put(1419.0,123.0){\rule[-0.200pt]{4.818pt}{0.400pt}}
\put(201.0,221.0){\rule[-0.200pt]{4.818pt}{0.400pt}}
\put(181,221){\makebox(0,0)[r]{ 0.02}}
\put(1419.0,221.0){\rule[-0.200pt]{4.818pt}{0.400pt}}
\put(201.0,320.0){\rule[-0.200pt]{4.818pt}{0.400pt}}
\put(181,320){\makebox(0,0)[r]{ 0.04}}
\put(1419.0,320.0){\rule[-0.200pt]{4.818pt}{0.400pt}}
\put(201.0,418.0){\rule[-0.200pt]{4.818pt}{0.400pt}}
\put(181,418){\makebox(0,0)[r]{ 0.06}}
\put(1419.0,418.0){\rule[-0.200pt]{4.818pt}{0.400pt}}
\put(201.0,516.0){\rule[-0.200pt]{4.818pt}{0.400pt}}
\put(181,516){\makebox(0,0)[r]{ 0.08}}
\put(1419.0,516.0){\rule[-0.200pt]{4.818pt}{0.400pt}}
\put(201.0,614.0){\rule[-0.200pt]{4.818pt}{0.400pt}}
\put(181,614){\makebox(0,0)[r]{ 0.1}}
\put(1419.0,614.0){\rule[-0.200pt]{4.818pt}{0.400pt}}
\put(201.0,713.0){\rule[-0.200pt]{4.818pt}{0.400pt}}
\put(181,713){\makebox(0,0)[r]{ 0.12}}
\put(1419.0,713.0){\rule[-0.200pt]{4.818pt}{0.400pt}}
\put(201.0,811.0){\rule[-0.200pt]{4.818pt}{0.400pt}}
\put(181,811){\makebox(0,0)[r]{ 0.14}}
\put(1419.0,811.0){\rule[-0.200pt]{4.818pt}{0.400pt}}
\put(201.0,123.0){\rule[-0.200pt]{0.400pt}{4.818pt}}
\put(201,82){\makebox(0,0){ 0.5}}
\put(201.0,840.0){\rule[-0.200pt]{0.400pt}{4.818pt}}
\put(449.0,123.0){\rule[-0.200pt]{0.400pt}{4.818pt}}
\put(449,82){\makebox(0,0){ 0.6}}
\put(449.0,840.0){\rule[-0.200pt]{0.400pt}{4.818pt}}
\put(696.0,123.0){\rule[-0.200pt]{0.400pt}{4.818pt}}
\put(696,82){\makebox(0,0){ 0.7}}
\put(696.0,840.0){\rule[-0.200pt]{0.400pt}{4.818pt}}
\put(944.0,123.0){\rule[-0.200pt]{0.400pt}{4.818pt}}
\put(944,82){\makebox(0,0){ 0.8}}
\put(944.0,840.0){\rule[-0.200pt]{0.400pt}{4.818pt}}
\put(1191.0,123.0){\rule[-0.200pt]{0.400pt}{4.818pt}}
\put(1191,82){\makebox(0,0){ 0.9}}
\put(1191.0,840.0){\rule[-0.200pt]{0.400pt}{4.818pt}}
\put(1439.0,123.0){\rule[-0.200pt]{0.400pt}{4.818pt}}
\put(1439,82){\makebox(0,0){ 1}}
\put(1439.0,840.0){\rule[-0.200pt]{0.400pt}{4.818pt}}
\put(201.0,123.0){\rule[-0.200pt]{298.234pt}{0.400pt}}
\put(1439.0,123.0){\rule[-0.200pt]{0.400pt}{177.543pt}}
\put(201.0,860.0){\rule[-0.200pt]{298.234pt}{0.400pt}}
\put(201.0,123.0){\rule[-0.200pt]{0.400pt}{177.543pt}}
\put(40,491){\makebox(0,0){$\sigma$}}
\put(820,21){\makebox(0,0){$R$}}
\put(1279,820){\makebox(0,0)[r]{Example 2 without rounding}}
\put(1299.0,820.0){\rule[-0.200pt]{24.090pt}{0.400pt}}
\put(201,427){\usebox{\plotpoint}}
\put(201.0,427.0){\rule[-0.200pt]{28.908pt}{0.400pt}}
\put(321.0,427.0){\rule[-0.200pt]{0.400pt}{104.310pt}}
\put(322,860){\usebox{\plotpoint}}
\put(322.0,577.0){\rule[-0.200pt]{0.400pt}{68.175pt}}
\put(322.0,577.0){\rule[-0.200pt]{61.429pt}{0.400pt}}
\put(577.0,423.0){\rule[-0.200pt]{0.400pt}{37.099pt}}
\put(577.0,423.0){\rule[-0.200pt]{29.872pt}{0.400pt}}
\put(701.0,423.0){\rule[-0.200pt]{0.400pt}{1.927pt}}
\put(701.0,431.0){\rule[-0.200pt]{28.426pt}{0.400pt}}
\put(819.0,431.0){\rule[-0.200pt]{0.400pt}{24.572pt}}
\put(819.0,533.0){\rule[-0.200pt]{29.631pt}{0.400pt}}
\put(942.0,533.0){\rule[-0.200pt]{0.400pt}{30.835pt}}
\put(942.0,661.0){\rule[-0.200pt]{61.670pt}{0.400pt}}
\put(1198.0,661.0){\rule[-0.200pt]{0.400pt}{22.885pt}}
\put(1198.0,756.0){\rule[-0.200pt]{58.057pt}{0.400pt}}
\put(201.0,123.0){\rule[-0.200pt]{298.234pt}{0.400pt}}
\put(1439.0,123.0){\rule[-0.200pt]{0.400pt}{177.543pt}}
\put(201.0,860.0){\rule[-0.200pt]{298.234pt}{0.400pt}}
\put(201.0,123.0){\rule[-0.200pt]{0.400pt}{177.543pt}}
\end{picture}
\end{center}
\caption{Example 2: Deviation from fair voting without rounding.}
\label{f_ex2_deviation4}
\end{figure}
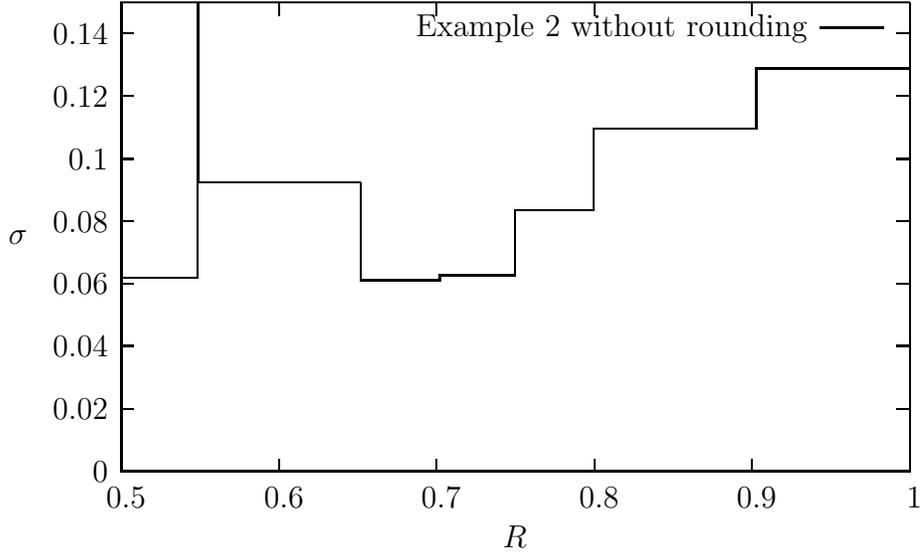

Again, the 4 digit voting weights are exact, and for these unrounded weights, we can list all possible coalitions
with their voting weights as shown in table~\ref{t_ex2_coalitions4}, and calculate $\sigma$ for all values of the
quota $R$ (table~\ref{t_ex2_deviation4} and figure~\ref{f_ex2_deviation4}).

\begin{table}[ht]
\begin{center}
\begin{tabular}{|c|r||c|r|}
\hline
Countries & $\sum v_i$ & Countries & $\sum v_i$ \\
\hline
1, 2, 3, 4 & 1.0000 & 1 & 0.4514 \\
1, 2, 3 & 0.9026 & 2, 3 & 0.4512 \\
1, 2, 4 & 0.7994 & 2, 4 & 0.3480 \\
1, 3, 4 & 0.7494 & 3, 4 & 0.2980 \\
1, 2 & 0.7020 & 2 & 0.2506 \\
1, 3 & 0.6520 & 3 & 0.2006 \\
1, 4 & 0.5488 & 4 & 0.0974 \\
2, 3, 4 & 0.5486 & None & 0.0000 \\
\hline
\end{tabular}
\end{center}
\caption{Example 2: Coalitions and their voting weights without rounding.}
\label{t_ex2_coalitions4}
\end{table}

\begin{table}[ht]
\begin{center}
\begin{tabular}{|r|r||c|c|c|c||r|}
\hline
\multicolumn{1}{|c|}{$R>$} & \multicolumn{1}{c||}{$R\leq$} & $\eta_1$ & $\eta_2$
 & $\eta_3$ & $\eta_4$ & \multicolumn{1}{c|}{$\sigma$} \\
\hline
 $0.9026$ & $1.0000$ & $1$ & $1$ & $1$ & $1$ & $0.128734$ \\
 $0.7994$ & $0.9026$ & $2$ & $2$ & $2$ & $0$ & $0.109418$ \\
 $0.7494$ & $0.7994$ & $3$ & $3$ & $1$ & $1$ & $0.083351$ \\
 $0.7020$ & $0.7494$ & $4$ & $2$ & $2$ & $2$ & $0.062709$ \\
 $0.6520$ & $0.7020$ & $5$ & $3$ & $1$ & $1$ & $0.061093$ \\
 $0.5488$ & $0.6520$ & $6$ & $2$ & $2$ & $0$ & $0.092371$ \\
 $0.5486$ & $0.5488$ & $7$ & $1$ & $1$ & $1$ & $0.153793$ \\
 $0.4514$ & $0.5486$ & $6$ & $2$ & $2$ & $2$ & $0.061960$ \\
 $0.4512$ & $0.4514$ & $7$ & $1$ & $1$ & $1$ & $0.153793$ \\
 $0.3480$ & $0.4512$ & $6$ & $2$ & $2$ & $0$ & $0.092371$ \\
 $0.2980$ & $0.3480$ & $5$ & $3$ & $1$ & $1$ & $0.061093$ \\
 $0.2506$ & $0.2980$ & $4$ & $2$ & $2$ & $2$ & $0.062709$ \\
 $0.2006$ & $0.2506$ & $3$ & $3$ & $1$ & $1$ & $0.083351$ \\
 $0.0974$ & $0.2006$ & $2$ & $2$ & $2$ & $0$ & $0.109418$ \\
 $0.0000$ & $0.0974$ & $1$ & $1$ & $1$ & $1$ & $0.128734$ \\
\hline
\end{tabular}
\end{center}
\caption{Example 2: Deviation from fair voting without rounding.}
\label{t_ex2_deviation4}
\end{table}

We can read off the minimum value of $\sigma$, which is
\begin{equation}
  \sigma_{\rm min}=0.061093,
\end{equation}
and occurs at
\begin{equation}
  0.6520<R\leq 0.7020.
\end{equation}

\begin{figure}[ht]
\begin{center}
% GNUPLOT: LaTeX picture
\setlength{\unitlength}{0.240900pt}
\ifx\plotpoint\undefined\newsavebox{\plotpoint}\fi
\sbox{\plotpoint}{\rule[-0.200pt]{0.400pt}{0.400pt}}%
\begin{picture}(1500,900)(0,0)
\sbox{\plotpoint}{\rule[-0.200pt]{0.400pt}{0.400pt}}%
\put(201.0,123.0){\rule[-0.200pt]{4.818pt}{0.400pt}}
\put(181,123){\makebox(0,0)[r]{ 0}}
\put(1419.0,123.0){\rule[-0.200pt]{4.818pt}{0.400pt}}
\put(201.0,221.0){\rule[-0.200pt]{4.818pt}{0.400pt}}
\put(181,221){\makebox(0,0)[r]{ 0.02}}
\put(1419.0,221.0){\rule[-0.200pt]{4.818pt}{0.400pt}}
\put(201.0,320.0){\rule[-0.200pt]{4.818pt}{0.400pt}}
\put(181,320){\makebox(0,0)[r]{ 0.04}}
\put(1419.0,320.0){\rule[-0.200pt]{4.818pt}{0.400pt}}
\put(201.0,418.0){\rule[-0.200pt]{4.818pt}{0.400pt}}
\put(181,418){\makebox(0,0)[r]{ 0.06}}
\put(1419.0,418.0){\rule[-0.200pt]{4.818pt}{0.400pt}}
\put(201.0,516.0){\rule[-0.200pt]{4.818pt}{0.400pt}}
\put(181,516){\makebox(0,0)[r]{ 0.08}}
\put(1419.0,516.0){\rule[-0.200pt]{4.818pt}{0.400pt}}
\put(201.0,614.0){\rule[-0.200pt]{4.818pt}{0.400pt}}
\put(181,614){\makebox(0,0)[r]{ 0.1}}
\put(1419.0,614.0){\rule[-0.200pt]{4.818pt}{0.400pt}}
\put(201.0,713.0){\rule[-0.200pt]{4.818pt}{0.400pt}}
\put(181,713){\makebox(0,0)[r]{ 0.12}}
\put(1419.0,713.0){\rule[-0.200pt]{4.818pt}{0.400pt}}
\put(201.0,811.0){\rule[-0.200pt]{4.818pt}{0.400pt}}
\put(181,811){\makebox(0,0)[r]{ 0.14}}
\put(1419.0,811.0){\rule[-0.200pt]{4.818pt}{0.400pt}}
\put(201.0,123.0){\rule[-0.200pt]{0.400pt}{4.818pt}}
\put(201,82){\makebox(0,0){ 0.5}}
\put(201.0,840.0){\rule[-0.200pt]{0.400pt}{4.818pt}}
\put(449.0,123.0){\rule[-0.200pt]{0.400pt}{4.818pt}}
\put(449,82){\makebox(0,0){ 0.6}}
\put(449.0,840.0){\rule[-0.200pt]{0.400pt}{4.818pt}}
\put(696.0,123.0){\rule[-0.200pt]{0.400pt}{4.818pt}}
\put(696,82){\makebox(0,0){ 0.7}}
\put(696.0,840.0){\rule[-0.200pt]{0.400pt}{4.818pt}}
\put(944.0,123.0){\rule[-0.200pt]{0.400pt}{4.818pt}}
\put(944,82){\makebox(0,0){ 0.8}}
\put(944.0,840.0){\rule[-0.200pt]{0.400pt}{4.818pt}}
\put(1191.0,123.0){\rule[-0.200pt]{0.400pt}{4.818pt}}
\put(1191,82){\makebox(0,0){ 0.9}}
\put(1191.0,840.0){\rule[-0.200pt]{0.400pt}{4.818pt}}
\put(1439.0,123.0){\rule[-0.200pt]{0.400pt}{4.818pt}}
\put(1439,82){\makebox(0,0){ 1}}
\put(1439.0,840.0){\rule[-0.200pt]{0.400pt}{4.818pt}}
\put(201.0,123.0){\rule[-0.200pt]{298.234pt}{0.400pt}}
\put(1439.0,123.0){\rule[-0.200pt]{0.400pt}{177.543pt}}
\put(201.0,860.0){\rule[-0.200pt]{298.234pt}{0.400pt}}
\put(201.0,123.0){\rule[-0.200pt]{0.400pt}{177.543pt}}
\put(40,491){\makebox(0,0){$\sigma$}}
\put(820,21){\makebox(0,0){$R$}}
\put(1279,820){\makebox(0,0)[r]{Example 2 with 3-digit rounding}}
\put(1299.0,820.0){\rule[-0.200pt]{24.090pt}{0.400pt}}
\put(201,427){\usebox{\plotpoint}}
\put(201.0,427.0){\rule[-0.200pt]{28.667pt}{0.400pt}}
\put(320.0,275.0){\rule[-0.200pt]{0.400pt}{36.617pt}}
\put(320.0,275.0){\rule[-0.200pt]{0.482pt}{0.400pt}}
\put(322.0,275.0){\rule[-0.200pt]{0.400pt}{72.752pt}}
\put(322.0,577.0){\rule[-0.200pt]{61.429pt}{0.400pt}}
\put(577.0,423.0){\rule[-0.200pt]{0.400pt}{37.099pt}}
\put(577.0,423.0){\rule[-0.200pt]{29.872pt}{0.400pt}}
\put(701.0,423.0){\rule[-0.200pt]{0.400pt}{1.927pt}}
\put(701.0,431.0){\rule[-0.200pt]{28.185pt}{0.400pt}}
\put(818.0,431.0){\rule[-0.200pt]{0.400pt}{24.572pt}}
\put(818.0,533.0){\rule[-0.200pt]{29.631pt}{0.400pt}}
\put(941.0,533.0){\rule[-0.200pt]{0.400pt}{30.835pt}}
\put(941.0,661.0){\rule[-0.200pt]{62.152pt}{0.400pt}}
\put(1199.0,661.0){\rule[-0.200pt]{0.400pt}{22.885pt}}
\put(1199.0,756.0){\rule[-0.200pt]{57.816pt}{0.400pt}}
\put(201.0,123.0){\rule[-0.200pt]{298.234pt}{0.400pt}}
\put(1439.0,123.0){\rule[-0.200pt]{0.400pt}{177.543pt}}
\put(201.0,860.0){\rule[-0.200pt]{298.234pt}{0.400pt}}
\put(201.0,123.0){\rule[-0.200pt]{0.400pt}{177.543pt}}
\end{picture}
\end{center}
\caption{Example 2: Deviation from fair voting with 3 digit rounding.}
\label{f_ex2_deviation3}
\end{figure}
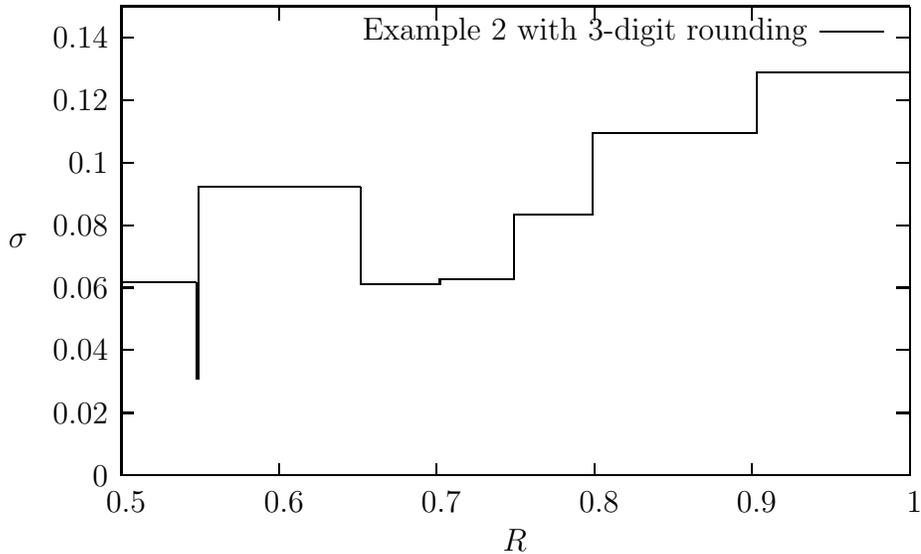

As in the previous example, the situation changes significantly if the voting weights are rounded to three digits.
Looking at the coalitions and their voting weights (table~\ref{t_ex2_coalitions3}), it can be seen that the
coalition consisting of countries 1 and 4, and the coalition consisting of countries 2, 3 and 4 have swapped places.
This is reflected in the function $\sigma$, as shown in table~\ref{t_ex2_deviation3} and figure~\ref{f_ex2_deviation3}.

\begin{table}[ht]
\begin{center}
\begin{tabular}{|c|r||c|r|}
\hline
Countries & $\sum v_i$ & Countries & $\sum v_i$ \\
\hline
1, 2, 3, 4 & 1.000 & 2, 3 & 0.452 \\
1, 2, 3 & 0.903 & 1 & 0.451 \\
1, 2, 4 & 0.799 & 2, 4 & 0.348 \\
1, 3, 4 & 0.749 & 3, 4 & 0.298 \\
1, 2 & 0.702 & 2 & 0.251 \\
1, 3 & 0.652 & 3 & 0.201 \\
2, 3, 4 & 0.549 & 4 & 0.097 \\
1, 4 & 0.548 & None & 0.000 \\
\hline
\end{tabular}
\end{center}
\caption{Example 2: Coalitions and their voting weights with 3 digit rounding.}
\label{t_ex2_coalitions3}
\end{table}

\begin{table}[ht]
\begin{center}
\begin{tabular}{|r|r||c|c|c|c||r|}
\hline
\multicolumn{1}{|c|}{$R>$} & \multicolumn{1}{c||}{$R\leq$} & $\eta_1$ & $\eta_2$ & $\eta_3$ & $\eta_4$ & \multicolumn{1}{c|}{$\sigma$} \\
\hline
 $0.903$ & $1.000$ & $1$ & $1$ & $1$ & $1$ & $0.128734$ \\
 $0.799$ & $0.903$ & $2$ & $2$ & $2$ & $0$ & $0.109418$ \\
 $0.749$ & $0.799$ & $3$ & $3$ & $1$ & $1$ & $0.083351$ \\
 $0.702$ & $0.749$ & $4$ & $2$ & $2$ & $2$ & $0.062709$ \\
 $0.652$ & $0.702$ & $5$ & $3$ & $1$ & $1$ & $0.061093$ \\
 $0.549$ & $0.652$ & $6$ & $2$ & $2$ & $0$ & $0.092371$ \\
 $0.548$ & $0.549$ & $5$ & $3$ & $3$ & $1$ & $0.031004$ \\
 $0.452$ & $0.548$ & $6$ & $2$ & $2$ & $2$ & $0.061960$ \\
 $0.451$ & $0.452$ & $5$ & $3$ & $3$ & $1$ & $0.031004$ \\
 $0.348$ & $0.451$ & $6$ & $2$ & $2$ & $0$ & $0.092371$ \\
 $0.298$ & $0.348$ & $5$ & $3$ & $1$ & $1$ & $0.061093$ \\
 $0.251$ & $0.298$ & $4$ & $2$ & $2$ & $2$ & $0.062709$ \\
 $0.201$ & $0.251$ & $3$ & $3$ & $1$ & $1$ & $0.083351$ \\
 $0.097$ & $0.201$ & $2$ & $2$ & $2$ & $0$ & $0.109418$ \\
 $0.000$ & $0.097$ & $1$ & $1$ & $1$ & $1$ & $0.128734$ \\
\hline
\end{tabular}
\end{center}
\caption{Example 2: Deviation from fair voting with 3 digit rounding.}
\label{t_ex2_deviation3}
\end{table}

Reading off the minimum value of $\sigma$, we now obtain
\begin{equation}
  \sigma_{\rm min}=0.031004,
\end{equation}
occurring at
\begin{equation}
  0.548<R\leq 0.549.
\end{equation}
In this case, the rounding has actually improved the approximation, the optimal value of $\sigma$ has almost halved. This means
that rounding can be used as a tool in finding the optimal quota, ie $\sigma$ can be calculated for several rounding schemes
(1 digit, 2 digits, \ldots), and the minimum value of $\sigma_{\rm min}$ found will then
determine the optimal rounding of the voting
weights, and the quota $R$ that should be chosen.

\section{Voting in the Council of the EU}
\label{s_euvoting}

\subsection{Current voting arrangements}
\label{ss_current}

The current voting procedures in the Council of the EU are laid down in the Treaty of Nice (\cite{NiceTreaty}) and the accession
treaties for the member states that joined the EU after the introduction of the Nice voting scheme (\cite{2003AccTreaty,%
2005AccTreaty}). Under this voting scheme,
each of the currently 27 member states is assigned a certain number of votes, and the votes for all member states add up to
345. For a proposal to be passed, the following criteria must be met:
\begin{enumerate}
\item The votes of the member states voting ``yes'' must add up to at least 255.
\item More than half of the member states (ie currently 14 countries) vote ``yes''. If the Council is not acting on a proposal
by the European Commission, at least two thirds of the member states (ie currently 18 countries) must vote ``yes''.
\item The countries voting ``yes'' must represent at least 62\% of the population of the EU.
\end{enumerate}
The last condition is only applied if a member state requests it, and there are only very few coalitions that satisfy the first two
but not the third condition.

The current population numbers (\cite{EuroStat}) and votes for each country are given in table~\ref{t_population}. At the time of writing,
the 2007 population data from Belgium, Ireland and Luxembourg were not available, and for these countries,
the 2006 data will be used in the calculations in this and the next two subsections.

\begin{center}
\begin{longtable}{lc@{}lc}
\hline
\multicolumn{1}{|c}{Country} & \multicolumn{2}{|c}{Population} & \multicolumn{1}{|c|}{Votes (Nice)}\\
\hline
\endfirsthead

%\multicolumn{4}{l}{\tablename\ \thetable{} -- continued from previous page}\\
%\hline
%\multicolumn{1}{|c}{Country} & \multicolumn{2}{|c}{Population} & \multicolumn{1}{|c|}{Votes (Nice)}\\
%\hline
%\endhead

\hline
\multicolumn{4}{r}{continued on next page}
\endfoot
 
\endlastfoot

\multicolumn{1}{|l}{Germany} & \multicolumn{1}{|r}{$82,310,995$} & $^b$ & \multicolumn{1}{|r|}{$29$} \\
\multicolumn{1}{|l}{France} & \multicolumn{1}{|r}{$63,392,140$} & $^b$ & \multicolumn{1}{|r|}{$29$} \\
\multicolumn{1}{|l}{UK} & \multicolumn{1}{|r}{$60,798,438$} & $^b$ & \multicolumn{1}{|r|}{$29$} \\
\multicolumn{1}{|l}{Italy} & \multicolumn{1}{|r}{$59,131,287$} & & \multicolumn{1}{|r|}{$29$} \\
\multicolumn{1}{|l}{Spain} & \multicolumn{1}{|r}{$44,474,631$} & & \multicolumn{1}{|r|}{$27$} \\
\multicolumn{1}{|l}{Poland} & \multicolumn{1}{|r}{$38,125,479$} & & \multicolumn{1}{|r|}{$27$} \\
\multicolumn{1}{|l}{Romania} & \multicolumn{1}{|r}{$21,565,119$} & $^b$ & \multicolumn{1}{|r|}{$14$} \\
\multicolumn{1}{|l}{Netherlands} & \multicolumn{1}{|r}{$16,357,992$} & & \multicolumn{1}{|r|}{$13$} \\
\multicolumn{1}{|l}{Greece} & \multicolumn{1}{|r}{$11,170,957$} & $^b$ & \multicolumn{1}{|r|}{$12$} \\
\multicolumn{1}{|l}{Portugal} & \multicolumn{1}{|r}{$10,599,095$} & & \multicolumn{1}{|r|}{$12$} \\
\multicolumn{1}{|l}{Belgium} & \multicolumn{1}{|r}{$10,511,382$} & $^a$ & \multicolumn{1}{|r|}{$12$} \\
\multicolumn{1}{|l}{Czech Republic} & \multicolumn{1}{|r}{$10,287,189$} & & \multicolumn{1}{|r|}{$12$} \\
\multicolumn{1}{|l}{Hungary} & \multicolumn{1}{|r}{$10,064,000$} & $^b$ & \multicolumn{1}{|r|}{$12$} \\
\multicolumn{1}{|l}{Sweden} & \multicolumn{1}{|r}{$9,113,257$} & & \multicolumn{1}{|r|}{$10$} \\
\multicolumn{1}{|l}{Austria} & \multicolumn{1}{|r}{$8,298,923$} & $^b$ & \multicolumn{1}{|r|}{$10$} \\
\multicolumn{1}{|l}{Bulgaria} & \multicolumn{1}{|r}{$7,679,290$} & & \multicolumn{1}{|r|}{$10$} \\
\multicolumn{1}{|l}{Denmark} & \multicolumn{1}{|r}{$5,447,084$} & & \multicolumn{1}{|r|}{$7$} \\
\multicolumn{1}{|l}{Slovakia} & \multicolumn{1}{|r}{$5,393,637$} & & \multicolumn{1}{|r|}{$7$} \\
\multicolumn{1}{|l}{Finland} & \multicolumn{1}{|r}{$5,276,955$} & & \multicolumn{1}{|r|}{$7$} \\
\multicolumn{1}{|l}{Ireland} & \multicolumn{1}{|r}{$4,209,019$} & $^{a,b}$ & \multicolumn{1}{|r|}{$7$} \\
\multicolumn{1}{|l}{Lithuania} & \multicolumn{1}{|r}{$3,384,879$} & & \multicolumn{1}{|r|}{$7$} \\
\multicolumn{1}{|l}{Latvia} & \multicolumn{1}{|r}{$2,281,305$} & & \multicolumn{1}{|r|}{$4$} \\
\multicolumn{1}{|l}{Slovenia} & \multicolumn{1}{|r}{$2,010,377$} & & \multicolumn{1}{|r|}{$4$} \\
\multicolumn{1}{|l}{Estonia} & \multicolumn{1}{|r}{$1,342,409$} & & \multicolumn{1}{|r|}{$4$} \\
\multicolumn{1}{|l}{Cyprus} & \multicolumn{1}{|r}{$778,537$} & $^b$ & \multicolumn{1}{|r|}{$4$} \\
\multicolumn{1}{|l}{Luxembourg} & \multicolumn{1}{|r}{$459,500$} & $^{a,b}$ & \multicolumn{1}{|r|}{$4$} \\
\multicolumn{1}{|l}{Malta} & \multicolumn{1}{|r}{$406,020$} & & \multicolumn{1}{|r|}{$3$} \\
\hline
\multicolumn{1}{|l}{Total} & \multicolumn{1}{|r}{$494,869,896$} & & \multicolumn{1}{|r|}{$345$} \\
\hline
$^a$ 2006 & & & \\
$^b$ provisional & & & \\
\caption{Population and votes under the Treaty of Nice.}
\label{t_population}
\end{longtable}
\end{center}

In 2004, a new voting scheme was suggested in the draft constitution for Europe (\cite{DraftConstitution}), which subsequently
failed to reach a majority in referenda in France and in the Netherlands. However, an Intergovernmental Conference in 2007 has
drafted a new Reform Treaty (\cite{ReformTreaty}), which, if accepted, will introduce the voting method suggested by the draft
constitution in November 2014. According to this voting scheme, a proposal is passed, if it satisfies the ``double majority''
conditions:
\begin{enumerate}
\item At least 55\% of the member states (currently 15) vote ``yes'', or, if the Council is not acting on a proposal by the
European Commission or the High Representative of the Union for Foreign Affairs and Security Policy, at least 72\% of the
member states (currently 20) vote ``yes''.
\item The countries voting ``yes'' represent at least 65\% of the population of the European Union.
\end{enumerate}
During a transition period (2014--2017), the Nice voting system can be used if this is requested by a member state.

The introduction of a square root voting system was suggested by Sweden at the Intergovernmental Conference 2000 (\cite{slomczynski2}),
and by Poland in 2007 (\cite{euobserver}), but was turned down by a majority of EU member states on both occasions. The scheme has
also been brought to the attention of governments in an open letter written by scientists, \cite{OpenLetter}.

\begin{figure}[ht]
\begin{center}
\includegraphics{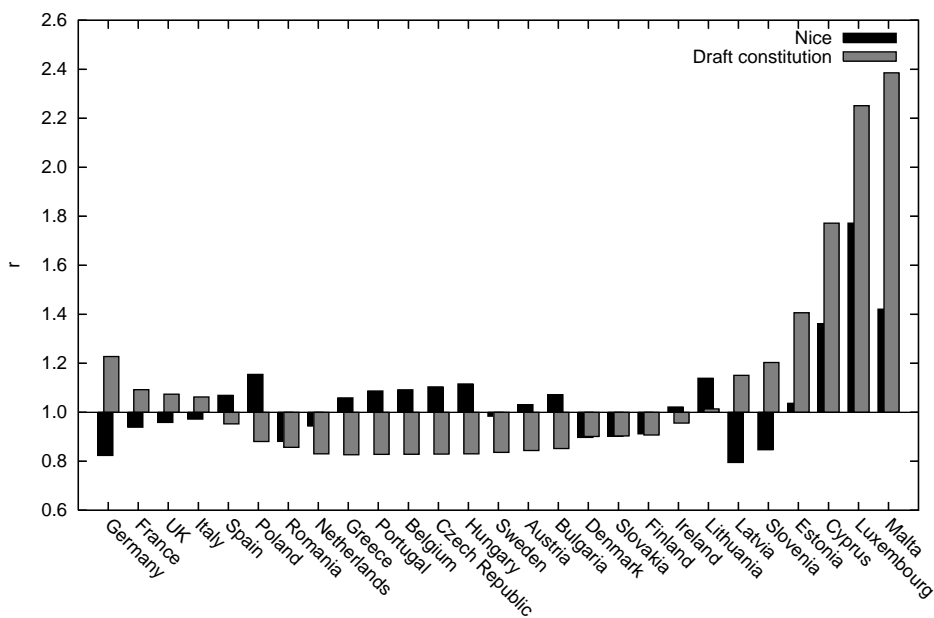}
\end{center}
\caption{Deviation from square root voting power for the Nice and draft Constitution voting schemes.}
\label{f_quotacomp}
\end{figure}

For both the Nice voting scheme and the draft constitution scheme, the relative Banzhaf indices $\beta_i$ can be calculated,
and the ratio
\begin{equation}
  \label{e_ratiodef}
  r = \frac{\beta_i}{\sqrt{N_i}/\sum_{j=1}^n \sqrt{N_j}}
\end{equation}
can be determined for each country, as a measure of ``fairness''. In a voting system that is fair in the Penrose sense,
\begin{equation}
  r = 1
\end{equation}
should hold for all countries. As shown in figure~\ref{f_quotacomp}, this is not the case: Both schemes favour smaller
countries (this is due to the allocation of a relatively large voting weight to small countries in the Nice scheme, and
the 55\% member quota in the draft constitution). Furthermore, the draft Constitution favours the largest countries due
to the 65\% population quota. This would be acceptable if these schemes were based on an alternative definition of
fairness and voting power, which does not seem to be the case.

For a detailed discussion of the shortcomings of the Nice and draft constitution schemes see \cite{Nurmi1}, \cite{Leech1}, \cite{Baldwin1}, \cite{Baldwin2},
\cite{Felsenthal1}, \cite{Felderer1}, \cite{Plechanovova1}, \cite{Barbera1}, \cite{Baldwin3}.

\subsection{Square root voting in the EU}
\label{ss_eusquareroot}

The square root voting scheme suggested by Poland follows the Jagiellonian Compromise paper \cite{slomczynski1}, which finds
(for the 27 member EU) an optimal quota of approximately $0.615$, in the sense that at this quota, the quantity $\sigma$
defined in section~\ref{s_squarerootvoting} is minimised. The purpose of this section is to find the ``exact'' minimum of
$\sigma$ using double precision numbers (ie 8 byte real numbers), and to investigate the effect of rounding on the minimum
found.

The method used to minimise $\sigma$ is the same as in section~\ref{s_rounding}. There are
\begin{equation}
  2^{27} = 134,217,728
\end{equation}
different coalitions in the 27 member EU. For each coalition, the square root voting weight is calculated, and the coalitions that
have a voting weight of $0.5$ or more are sorted into three subintervals of the interval $[0.5,1]$, and within each subinterval
are sorted according to their voting weights using the quicksort algorithm \cite{hoare1},\cite{hoare2}, in the version presented
in \cite{NumericalRecipes}. As in the four-country examples,
the relative Banzhaf indices are then calculated for all values of the quota $R$, and $\sigma$ is then determined as a function of
$R$. On a Pentium 4 PC, this procedure takes a few minutes. Searching through the values of $\sigma$, the minimum
is found for
\begin{equation}
  0.614966335781001 < R \leq 0.614966337885155,
\end{equation}
with a value of
\begin{equation}
  \sigma_{\rm min} = 4.334790\times 10^{-5}.
\end{equation}
The quantity $\sigma$ as a function of the quota $R$ around the minimum is shown in figures~\ref{f_deviation0_1}
to~\ref{f_deviation0_4}, with increasingly improved resolution, where resolution means the horizontal distance between
individual points.

\begin{figure}[ht]
\begin{center}
\input{deviation1.tex}
\end{center}
\caption{Square root voting in the EU: Deviation from fair voting.}
\label{f_deviation0_1}
\end{figure}

\begin{figure}[ht]
\begin{center}
\input{deviation2.tex}
\end{center}
\caption{Square root voting in the EU: Deviation from fair voting.}
\label{f_deviation0_2}
\end{figure}

\begin{figure}[ht]
\begin{center}
\input{deviation3.tex}
\end{center}
\caption{Square root voting in the EU: Deviation from fair voting.}
\label{f_deviation0_3}
\end{figure}

\begin{figure}[ht]
\begin{center}
% GNUPLOT: LaTeX picture
\setlength{\unitlength}{0.240900pt}
\ifx\plotpoint\undefined\newsavebox{\plotpoint}\fi
\sbox{\plotpoint}{\rule[-0.200pt]{0.400pt}{0.400pt}}%
\begin{picture}(1500,900)(0,0)
\sbox{\plotpoint}{\rule[-0.200pt]{0.400pt}{0.400pt}}%
\put(301.0,123.0){\rule[-0.200pt]{4.818pt}{0.400pt}}
\put(281,123){\makebox(0,0)[r]{ 4.334e-05}}
\put(1419.0,123.0){\rule[-0.200pt]{4.818pt}{0.400pt}}
\put(301.0,492.0){\rule[-0.200pt]{4.818pt}{0.400pt}}
\put(281,492){\makebox(0,0)[r]{ 4.335e-05}}
\put(1419.0,492.0){\rule[-0.200pt]{4.818pt}{0.400pt}}
\put(301.0,860.0){\rule[-0.200pt]{4.818pt}{0.400pt}}
\put(281,860){\makebox(0,0)[r]{ 4.336e-05}}
\put(1419.0,860.0){\rule[-0.200pt]{4.818pt}{0.400pt}}
\put(301.0,123.0){\rule[-0.200pt]{0.400pt}{4.818pt}}
\put(301,82){\makebox(0,0){0.61496632}}
\put(301.0,840.0){\rule[-0.200pt]{0.400pt}{4.818pt}}
\put(680.0,123.0){\rule[-0.200pt]{0.400pt}{4.818pt}}
\put(680,82){\makebox(0,0){0.61496633}}
\put(680.0,840.0){\rule[-0.200pt]{0.400pt}{4.818pt}}
\put(1060.0,123.0){\rule[-0.200pt]{0.400pt}{4.818pt}}
\put(1060,82){\makebox(0,0){0.61496634}}
\put(1060.0,840.0){\rule[-0.200pt]{0.400pt}{4.818pt}}
\put(1439.0,123.0){\rule[-0.200pt]{0.400pt}{4.818pt}}
\put(1439,82){\makebox(0,0){0.61496635}}
\put(1439.0,840.0){\rule[-0.200pt]{0.400pt}{4.818pt}}
\put(301.0,123.0){\rule[-0.200pt]{274.144pt}{0.400pt}}
\put(1439.0,123.0){\rule[-0.200pt]{0.400pt}{177.543pt}}
\put(301.0,860.0){\rule[-0.200pt]{274.144pt}{0.400pt}}
\put(301.0,123.0){\rule[-0.200pt]{0.400pt}{177.543pt}}
\put(40,491){\makebox(0,0){$\sigma$}}
\put(870,21){\makebox(0,0){$R$}}
\put(1279,820){\makebox(0,0)[r]{Exact function}}
\put(1299.0,820.0){\rule[-0.200pt]{24.090pt}{0.400pt}}
\put(1204.0,438.0){\rule[-0.200pt]{56.611pt}{0.400pt}}
\put(1204.0,438.0){\rule[-0.200pt]{0.400pt}{8.191pt}}
\put(1128.0,472.0){\rule[-0.200pt]{18.308pt}{0.400pt}}
\put(1128.0,472.0){\rule[-0.200pt]{0.400pt}{5.059pt}}
\put(1122.0,493.0){\rule[-0.200pt]{1.445pt}{0.400pt}}
\put(1122.0,437.0){\rule[-0.200pt]{0.400pt}{13.490pt}}
\put(979.0,437.0){\rule[-0.200pt]{34.449pt}{0.400pt}}
\put(979.0,414.0){\rule[-0.200pt]{0.400pt}{5.541pt}}
\put(900.0,414.0){\rule[-0.200pt]{19.031pt}{0.400pt}}
\put(900.0,414.0){\rule[-0.200pt]{0.400pt}{3.132pt}}
\put(653.0,427.0){\rule[-0.200pt]{59.502pt}{0.400pt}}
\put(653.0,427.0){\rule[-0.200pt]{0.400pt}{17.827pt}}
\put(534.0,501.0){\rule[-0.200pt]{28.667pt}{0.400pt}}
\put(534.0,501.0){\rule[-0.200pt]{0.400pt}{12.768pt}}
\put(516.0,554.0){\rule[-0.200pt]{4.336pt}{0.400pt}}
\put(516.0,524.0){\rule[-0.200pt]{0.400pt}{7.227pt}}
\put(409.0,524.0){\rule[-0.200pt]{25.776pt}{0.400pt}}
\put(409.0,517.0){\rule[-0.200pt]{0.400pt}{1.686pt}}
\put(384.0,517.0){\rule[-0.200pt]{6.022pt}{0.400pt}}
\put(384.0,517.0){\rule[-0.200pt]{0.400pt}{3.373pt}}
\put(301.0,531.0){\rule[-0.200pt]{19.995pt}{0.400pt}}
\put(301.0,123.0){\rule[-0.200pt]{274.144pt}{0.400pt}}
\put(1439.0,123.0){\rule[-0.200pt]{0.400pt}{177.543pt}}
\put(301.0,860.0){\rule[-0.200pt]{274.144pt}{0.400pt}}
\put(301.0,123.0){\rule[-0.200pt]{0.400pt}{177.543pt}}
\end{picture}
\end{center}
\caption{Square root voting in the EU: Deviation from fair voting, exact function.}
\label{f_deviation0_4}
\end{figure}

Only very fine resolution reveals the piecewise constant nature of the function $\sigma$, the
minimum is obtained in an interval of a length of just over $2\times 10^{-9}$. This means that rounding to less than
10 digits may
have a significant impact on the minimal value of $\sigma$ as well as on the range of quotas for which it is obtained.

\begin{figure}[ht]
\begin{center}
% GNUPLOT: LaTeX picture
\setlength{\unitlength}{0.240900pt}
\ifx\plotpoint\undefined\newsavebox{\plotpoint}\fi
\sbox{\plotpoint}{\rule[-0.200pt]{0.400pt}{0.400pt}}%
\begin{picture}(1500,900)(0,0)
\sbox{\plotpoint}{\rule[-0.200pt]{0.400pt}{0.400pt}}%
\put(241.0,123.0){\rule[-0.200pt]{4.818pt}{0.400pt}}
\put(221,123){\makebox(0,0)[r]{ 0}}
\put(1419.0,123.0){\rule[-0.200pt]{4.818pt}{0.400pt}}
\put(241.0,270.0){\rule[-0.200pt]{4.818pt}{0.400pt}}
\put(221,270){\makebox(0,0)[r]{ 0.0001}}
\put(1419.0,270.0){\rule[-0.200pt]{4.818pt}{0.400pt}}
\put(241.0,418.0){\rule[-0.200pt]{4.818pt}{0.400pt}}
\put(221,418){\makebox(0,0)[r]{ 0.0002}}
\put(1419.0,418.0){\rule[-0.200pt]{4.818pt}{0.400pt}}
\put(241.0,565.0){\rule[-0.200pt]{4.818pt}{0.400pt}}
\put(221,565){\makebox(0,0)[r]{ 0.0003}}
\put(1419.0,565.0){\rule[-0.200pt]{4.818pt}{0.400pt}}
\put(241.0,713.0){\rule[-0.200pt]{4.818pt}{0.400pt}}
\put(221,713){\makebox(0,0)[r]{ 0.0004}}
\put(1419.0,713.0){\rule[-0.200pt]{4.818pt}{0.400pt}}
\put(241.0,860.0){\rule[-0.200pt]{4.818pt}{0.400pt}}
\put(221,860){\makebox(0,0)[r]{ 0.0005}}
\put(1419.0,860.0){\rule[-0.200pt]{4.818pt}{0.400pt}}
\put(241.0,123.0){\rule[-0.200pt]{0.400pt}{4.818pt}}
\put(241,82){\makebox(0,0){ 2}}
\put(241.0,840.0){\rule[-0.200pt]{0.400pt}{4.818pt}}
\put(481.0,123.0){\rule[-0.200pt]{0.400pt}{4.818pt}}
\put(481,82){\makebox(0,0){ 4}}
\put(481.0,840.0){\rule[-0.200pt]{0.400pt}{4.818pt}}
\put(720.0,123.0){\rule[-0.200pt]{0.400pt}{4.818pt}}
\put(720,82){\makebox(0,0){ 6}}
\put(720.0,840.0){\rule[-0.200pt]{0.400pt}{4.818pt}}
\put(960.0,123.0){\rule[-0.200pt]{0.400pt}{4.818pt}}
\put(960,82){\makebox(0,0){ 8}}
\put(960.0,840.0){\rule[-0.200pt]{0.400pt}{4.818pt}}
\put(1199.0,123.0){\rule[-0.200pt]{0.400pt}{4.818pt}}
\put(1199,82){\makebox(0,0){ 10}}
\put(1199.0,840.0){\rule[-0.200pt]{0.400pt}{4.818pt}}
\put(1439.0,123.0){\rule[-0.200pt]{0.400pt}{4.818pt}}
\put(1439,82){\makebox(0,0){ 12}}
\put(1439.0,840.0){\rule[-0.200pt]{0.400pt}{4.818pt}}
\put(241.0,123.0){\rule[-0.200pt]{288.598pt}{0.400pt}}
\put(1439.0,123.0){\rule[-0.200pt]{0.400pt}{177.543pt}}
\put(241.0,860.0){\rule[-0.200pt]{288.598pt}{0.400pt}}
\put(241.0,123.0){\rule[-0.200pt]{0.400pt}{177.543pt}}
\put(40,491){\makebox(0,0){$\sigma_{\rm min}$}}
\put(840,21){\makebox(0,0){digit rounding}}
\put(361,464){\circle*{12}}
\put(481,204){\circle*{12}}
\put(600,188){\circle*{12}}
\put(720,187){\circle*{12}}
\put(840,187){\circle*{12}}
\put(960,187){\circle*{12}}
\put(1080,187){\circle*{12}}
\put(1199,187){\circle*{12}}
\put(1319,187){\circle*{12}}
\put(1439,187){\circle*{12}}
\put(241.0,123.0){\rule[-0.200pt]{288.598pt}{0.400pt}}
\put(1439.0,123.0){\rule[-0.200pt]{0.400pt}{177.543pt}}
\put(241.0,860.0){\rule[-0.200pt]{288.598pt}{0.400pt}}
\put(241.0,123.0){\rule[-0.200pt]{0.400pt}{177.543pt}}
\put(241.0,123.0){\rule[-0.200pt]{4.818pt}{0.400pt}}
\put(221,123){\makebox(0,0)[r]{ 0}}
\put(1419.0,123.0){\rule[-0.200pt]{4.818pt}{0.400pt}}
\put(241.0,270.0){\rule[-0.200pt]{4.818pt}{0.400pt}}
\put(221,270){\makebox(0,0)[r]{ 0.0001}}
\put(1419.0,270.0){\rule[-0.200pt]{4.818pt}{0.400pt}}
\put(241.0,418.0){\rule[-0.200pt]{4.818pt}{0.400pt}}
\put(221,418){\makebox(0,0)[r]{ 0.0002}}
\put(1419.0,418.0){\rule[-0.200pt]{4.818pt}{0.400pt}}
\put(241.0,565.0){\rule[-0.200pt]{4.818pt}{0.400pt}}
\put(221,565){\makebox(0,0)[r]{ 0.0003}}
\put(1419.0,565.0){\rule[-0.200pt]{4.818pt}{0.400pt}}
\put(241.0,713.0){\rule[-0.200pt]{4.818pt}{0.400pt}}
\put(221,713){\makebox(0,0)[r]{ 0.0004}}
\put(1419.0,713.0){\rule[-0.200pt]{4.818pt}{0.400pt}}
\put(241.0,860.0){\rule[-0.200pt]{4.818pt}{0.400pt}}
\put(221,860){\makebox(0,0)[r]{ 0.0005}}
\put(1419.0,860.0){\rule[-0.200pt]{4.818pt}{0.400pt}}
\put(241.0,123.0){\rule[-0.200pt]{0.400pt}{4.818pt}}
\put(241,82){\makebox(0,0){ 2}}
\put(241.0,840.0){\rule[-0.200pt]{0.400pt}{4.818pt}}
\put(481.0,123.0){\rule[-0.200pt]{0.400pt}{4.818pt}}
\put(481,82){\makebox(0,0){ 4}}
\put(481.0,840.0){\rule[-0.200pt]{0.400pt}{4.818pt}}
\put(720.0,123.0){\rule[-0.200pt]{0.400pt}{4.818pt}}
\put(720,82){\makebox(0,0){ 6}}
\put(720.0,840.0){\rule[-0.200pt]{0.400pt}{4.818pt}}
\put(960.0,123.0){\rule[-0.200pt]{0.400pt}{4.818pt}}
\put(960,82){\makebox(0,0){ 8}}
\put(960.0,840.0){\rule[-0.200pt]{0.400pt}{4.818pt}}
\put(1199.0,123.0){\rule[-0.200pt]{0.400pt}{4.818pt}}
\put(1199,82){\makebox(0,0){ 10}}
\put(1199.0,840.0){\rule[-0.200pt]{0.400pt}{4.818pt}}
\put(1439.0,123.0){\rule[-0.200pt]{0.400pt}{4.818pt}}
\put(1439,82){\makebox(0,0){ 12}}
\put(1439.0,840.0){\rule[-0.200pt]{0.400pt}{4.818pt}}
\put(241.0,123.0){\rule[-0.200pt]{288.598pt}{0.400pt}}
\put(1439.0,123.0){\rule[-0.200pt]{0.400pt}{177.543pt}}
\put(241.0,860.0){\rule[-0.200pt]{288.598pt}{0.400pt}}
\put(241.0,123.0){\rule[-0.200pt]{0.400pt}{177.543pt}}
\put(40,491){\makebox(0,0){$\sigma_{\rm min}$}}
\put(840,21){\makebox(0,0){digit rounding}}
\put(241,187){\usebox{\plotpoint}}
\put(241.00,187.00){\usebox{\plotpoint}}
\put(261.76,187.00){\usebox{\plotpoint}}
\put(282.51,187.00){\usebox{\plotpoint}}
\put(303.27,187.00){\usebox{\plotpoint}}
\put(324.02,187.00){\usebox{\plotpoint}}
\put(344.78,187.00){\usebox{\plotpoint}}
\put(365.53,187.00){\usebox{\plotpoint}}
\put(386.29,187.00){\usebox{\plotpoint}}
\put(407.04,187.00){\usebox{\plotpoint}}
\put(427.80,187.00){\usebox{\plotpoint}}
\put(448.55,187.00){\usebox{\plotpoint}}
\put(469.31,187.00){\usebox{\plotpoint}}
\put(490.07,187.00){\usebox{\plotpoint}}
\put(510.82,187.00){\usebox{\plotpoint}}
\put(531.58,187.00){\usebox{\plotpoint}}
\put(552.33,187.00){\usebox{\plotpoint}}
\put(573.09,187.00){\usebox{\plotpoint}}
\put(593.84,187.00){\usebox{\plotpoint}}
\put(614.60,187.00){\usebox{\plotpoint}}
\put(635.35,187.00){\usebox{\plotpoint}}
\put(656.11,187.00){\usebox{\plotpoint}}
\put(676.87,187.00){\usebox{\plotpoint}}
\put(697.62,187.00){\usebox{\plotpoint}}
\put(718.38,187.00){\usebox{\plotpoint}}
\put(739.13,187.00){\usebox{\plotpoint}}
\put(759.89,187.00){\usebox{\plotpoint}}
\put(780.64,187.00){\usebox{\plotpoint}}
\put(801.40,187.00){\usebox{\plotpoint}}
\put(822.15,187.00){\usebox{\plotpoint}}
\put(842.91,187.00){\usebox{\plotpoint}}
\put(863.66,187.00){\usebox{\plotpoint}}
\put(884.42,187.00){\usebox{\plotpoint}}
\put(905.18,187.00){\usebox{\plotpoint}}
\put(925.93,187.00){\usebox{\plotpoint}}
\put(946.69,187.00){\usebox{\plotpoint}}
\put(967.44,187.00){\usebox{\plotpoint}}
\put(988.20,187.00){\usebox{\plotpoint}}
\put(1008.95,187.00){\usebox{\plotpoint}}
\put(1029.71,187.00){\usebox{\plotpoint}}
\put(1050.46,187.00){\usebox{\plotpoint}}
\put(1071.22,187.00){\usebox{\plotpoint}}
\put(1091.98,187.00){\usebox{\plotpoint}}
\put(1112.73,187.00){\usebox{\plotpoint}}
\put(1133.49,187.00){\usebox{\plotpoint}}
\put(1154.24,187.00){\usebox{\plotpoint}}
\put(1175.00,187.00){\usebox{\plotpoint}}
\put(1195.75,187.00){\usebox{\plotpoint}}
\put(1216.51,187.00){\usebox{\plotpoint}}
\put(1237.26,187.00){\usebox{\plotpoint}}
\put(1258.02,187.00){\usebox{\plotpoint}}
\put(1278.77,187.00){\usebox{\plotpoint}}
\put(1299.53,187.00){\usebox{\plotpoint}}
\put(1320.29,187.00){\usebox{\plotpoint}}
\put(1341.04,187.00){\usebox{\plotpoint}}
\put(1361.80,187.00){\usebox{\plotpoint}}
\put(1382.55,187.00){\usebox{\plotpoint}}
\put(1403.31,187.00){\usebox{\plotpoint}}
\put(1424.06,187.00){\usebox{\plotpoint}}
\put(1439,187){\usebox{\plotpoint}}
\put(241.0,123.0){\rule[-0.200pt]{288.598pt}{0.400pt}}
\put(1439.0,123.0){\rule[-0.200pt]{0.400pt}{177.543pt}}
\put(241.0,860.0){\rule[-0.200pt]{288.598pt}{0.400pt}}
\put(241.0,123.0){\rule[-0.200pt]{0.400pt}{177.543pt}}
\end{picture}
\end{center}
\caption{Minimum value of $\sigma$ from rounded voting weights. The dotted line is the value of $\sigma_{\rm min}$ calculated
from double precision weights.}
\label{f_deviation0_r}
\end{figure}
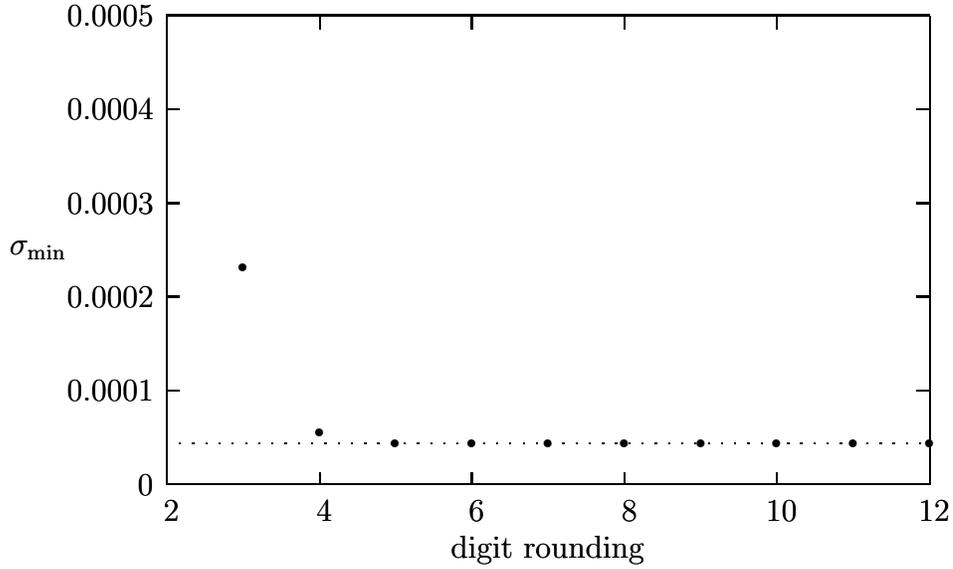

As pointed out in section~\ref{s_rounding}, the rounding effect may increase or decrease the deviation from fair
square root voting. Thus it is sensible to determine $\sigma_{\rm min}$ as described above for weights rounded
to $k$ digits,
$k=1,\ldots,k_{\rm max}\approx 10$, and choose the rounding scheme that returns the smallest value of $\sigma_{\min}$. This calculation
yields the results shown in figure~\ref{f_deviation0_r}.

Figure~\ref{f_deviation0_r} shows that the rounding effect on $\sigma_{\rm min}$ is very small provided the voting weights
are rounded to more than four decimal places, but the effect is clearly visible for the 4-digit rounding used
in the Jagiellonian Compromise paper \cite{slomczynski1}. Careful analysis shows that the minimum of $\sigma_{\rm min}$ is obtained with 7-digit rounding.
In this case, we obtain
\begin{equation}
  \sigma_{\rm min} = 4.334644\times 10^{-5},
\end{equation}
with an optimal quota $R$ in the interval
\begin{equation}
  0.6149670 < R \leq 0.6149671.
\end{equation}
Both the optimal quota and optimal rounding can change even with small changes in the population numbers. The same method
applied to the 2006 population data (\cite{EuroStat}) reveals that with those data, 4-digit rounding is optimal. This means
that the rounding scheme and the quota have to be adjusted annually, to make sure that the smallest possible deviation
from fair square root voting is obtained.

\begin{figure}[ht]
\begin{center}
\includegraphics{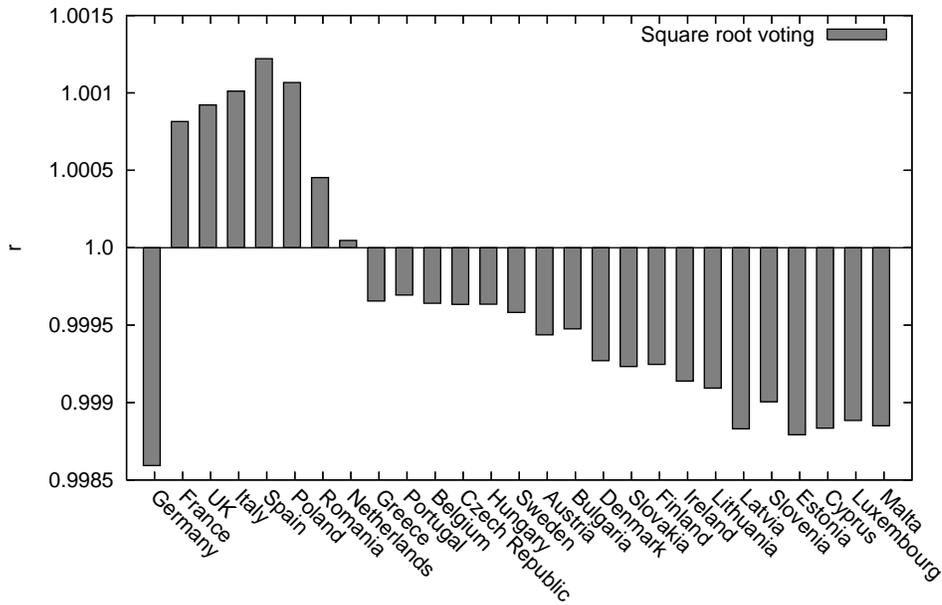}
\end{center}
\caption{Deviation from fair voting for all EU member states, for square root voting with 7-digit rounding.}
\label{f_quota0}
\end{figure}

Using 7-digit rounding on the 2007 data, the quotient $r$ defined in equation~\ref{e_ratiodef} can be calculated for
all 27 member states. This is shown in figure~\ref{f_quota0}.

A final quantity that is of interest is the efficiency $\epsilon$ of the voting scheme just described. It is
\begin{equation}
  \epsilon = 0.1644.
\end{equation}
This makes the square root method considerably more efficient than the Nice scheme ($\epsilon=0.02026$) and the
draft constitution ($\epsilon=0.1284$).

\subsection{Square root voting with a member quota}
\label{ss_eumemberquota}

Slomczy{\'n}ski and Zyczkowski point out (\cite{slomczynski1}) that it is possible to combine the square root voting
scheme with a requirement that a majority of member states has to vote ``yes'' for a proposal to be successful. In
fact, it is possible to include any member quota (14 countries, 15 countries, \ldots) in the square root method, which
could be used to make sure that the interests of smaller member states receive adequate attention. Member quotas can
also be employed as a threshold when voting on proposals not emanating from the European Commission, as is currently
the case with the 18 member quota under the Treaty of Nice.

\begin{figure}[ht]
\begin{center}
\input{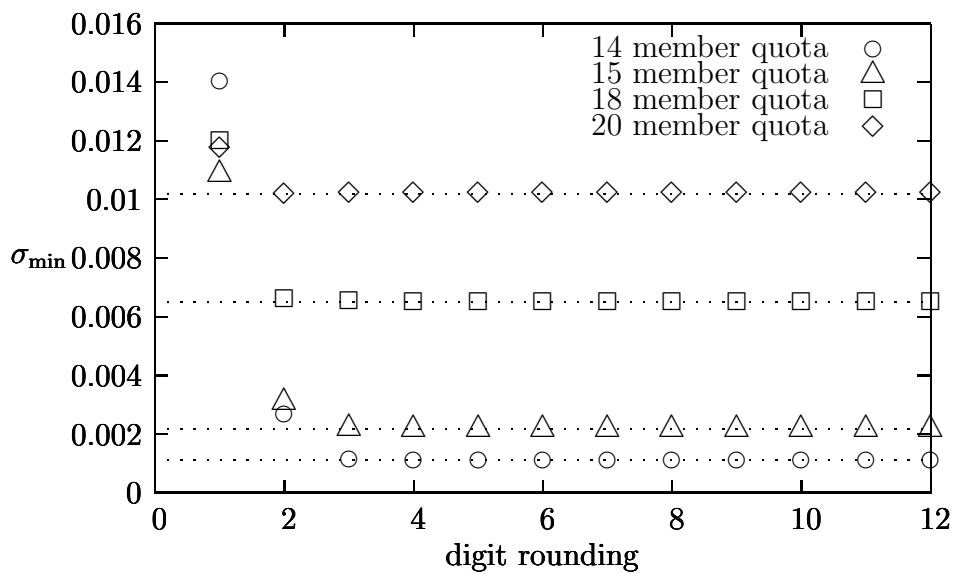}
\end{center}
\caption{Minimum value of $\sigma$ from rounded voting weights for different member quotas.
The dotted lines are the values of $\sigma_{\rm min}$ calculated
from double precision weights.}
\label{f_deviation14}
\end{figure}

As in the case without member quota, all possible coalitions are ordered according to their voting weights, and the
relative Banzhaf indices can then be calculated for all quotas $R$, simply by demanding that a coalition is successful
if its voting weight is at least equal to $R$, and the number of countries included in the coalition is at least
the member quota. The minimum deviation $\sigma_{\rm min}$ can then be found as above, and the calculation can be
repeated for 1-digit rounding, 2-digit rounding, etc. The results of this calculation for a number of member quotas
(14, 15, 18 and 20, to reflect the different member quotas in the Treaty of Nice and the draft constitution) are shown
in figure~\ref{f_deviation14}.

The data show that the optimal rounding for a member quota of 14, 15 or 18 is six digits, while for the
20 member quota, two digit rounding gives the smallest value of $\sigma_{\rm min}$. Table~\ref{t_memberquotasigma}
shows the optimal values of $\sigma_{\rm min}$, and the quota intervals where they occur.

\begin{table}[ht]
\begin{center}
\begin{tabular}{|r|r|r|r|}
\hline
Member quota & \multicolumn{1}{|c|}{$R>$} & \multicolumn{1}{c|}{$R\leq$} & \multicolumn{1}{c|}{$\sigma_{\rm min}$} \\
\hline
$14$ & $0.646660$ & $0.646661$ & $1.10654\times 10^{-3}$ \\
$15$ & $0.682884$ & $0.682885$ & $2.17490\times 10^{-3}$ \\
$18$ & $0.784222$ & $0.784223$ & $6.49661\times 10^{-3}$ \\
$20$ & $0.80$ & $0.81$ & $1.01626\times 10^{-2}$ \\
\hline
\end{tabular}
\end{center}
\caption{Optimal quota intervals and deviation from fair voting for different member quotas.}
\label{t_memberquotasigma}
\end{table}

\begin{figure}[ht]
\begin{center}
\includegraphics{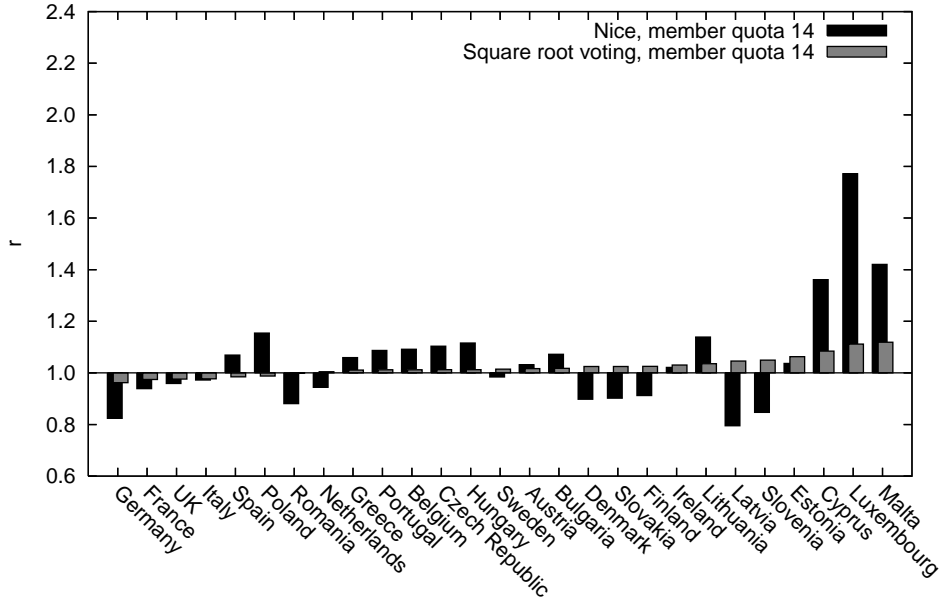}
\end{center}
\caption{The ratio $r$ for square root voting and Nice with a member quota of 14.}
\label{f_quota14}
\end{figure}

\begin{figure}[ht]
\begin{center}
\includegraphics{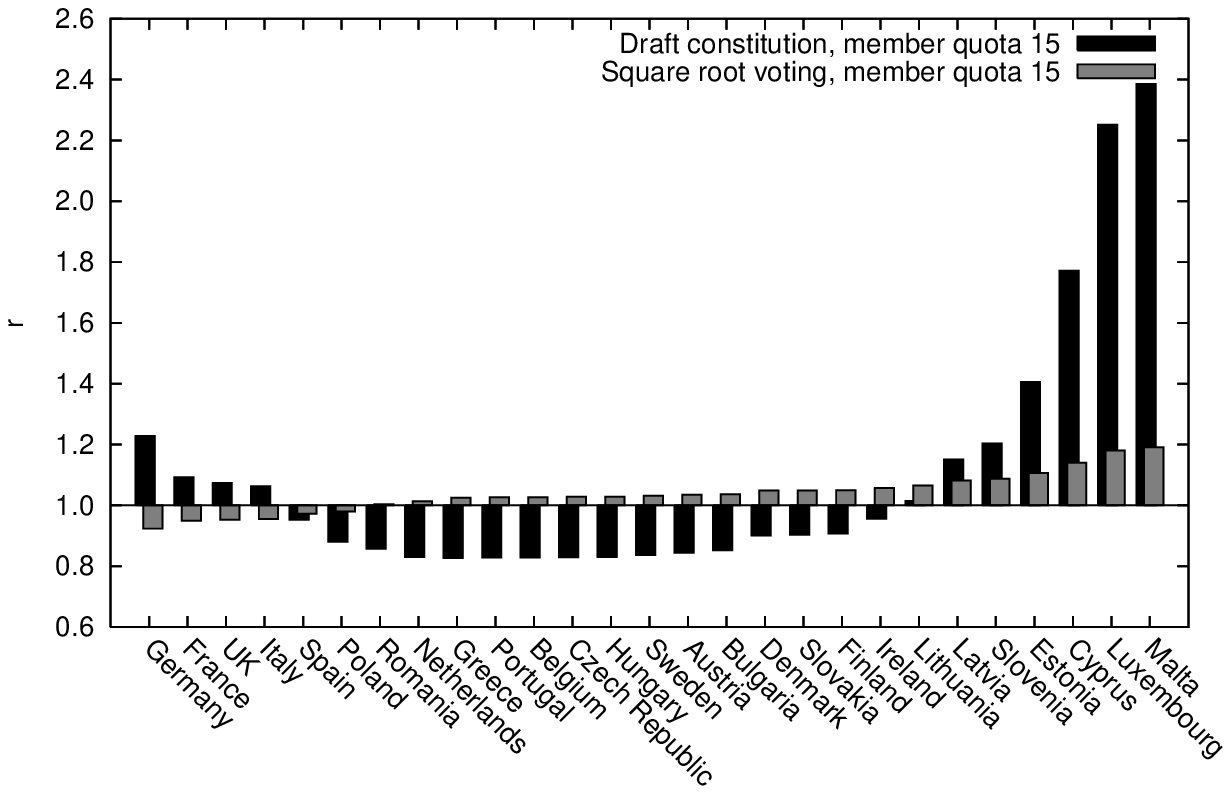}
\end{center}
\caption{The ratio $r$ for square root voting and the draft constitution with a member quota of 15.}
\label{f_quota15}
\end{figure}

\begin{figure}[ht]
\begin{center}
\includegraphics{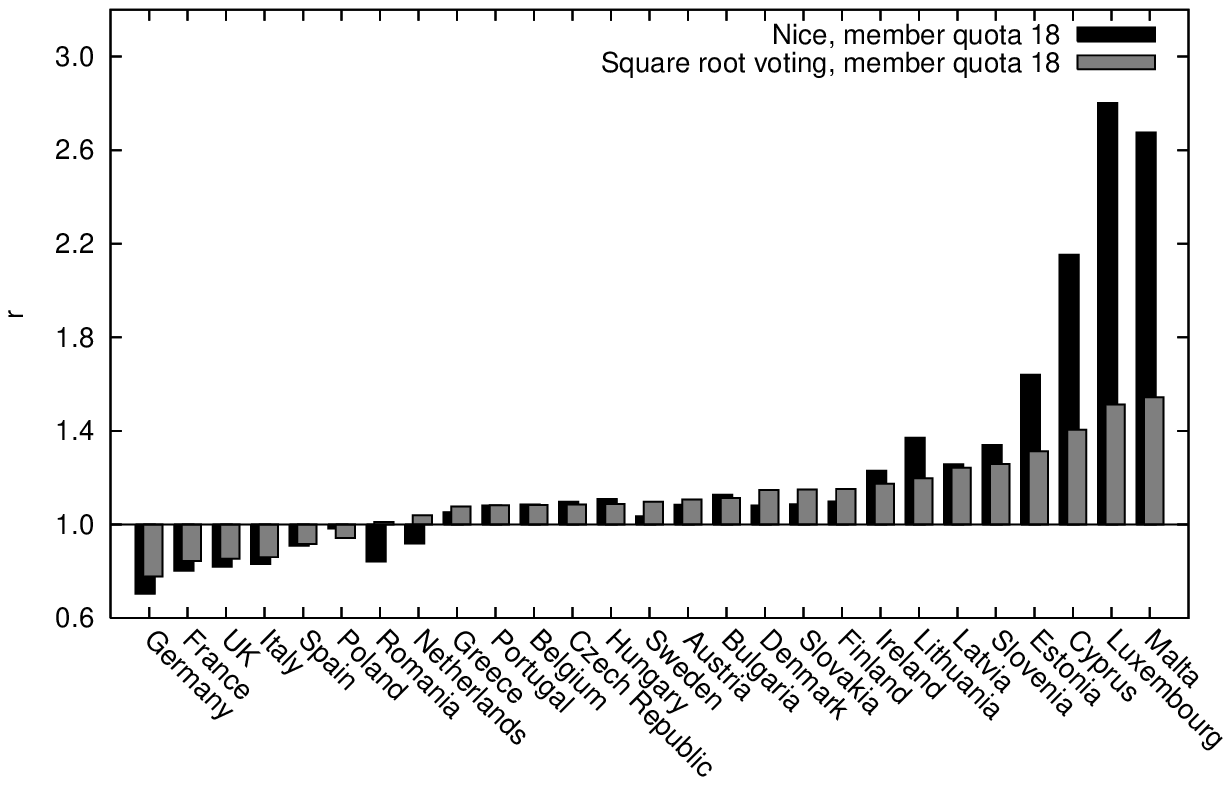}
\end{center}
\caption{The ratio $r$ for square root voting and Nice with a member quota of 18.}
\label{f_quota18}
\end{figure}

\begin{figure}[ht]
\begin{center}
\includegraphics{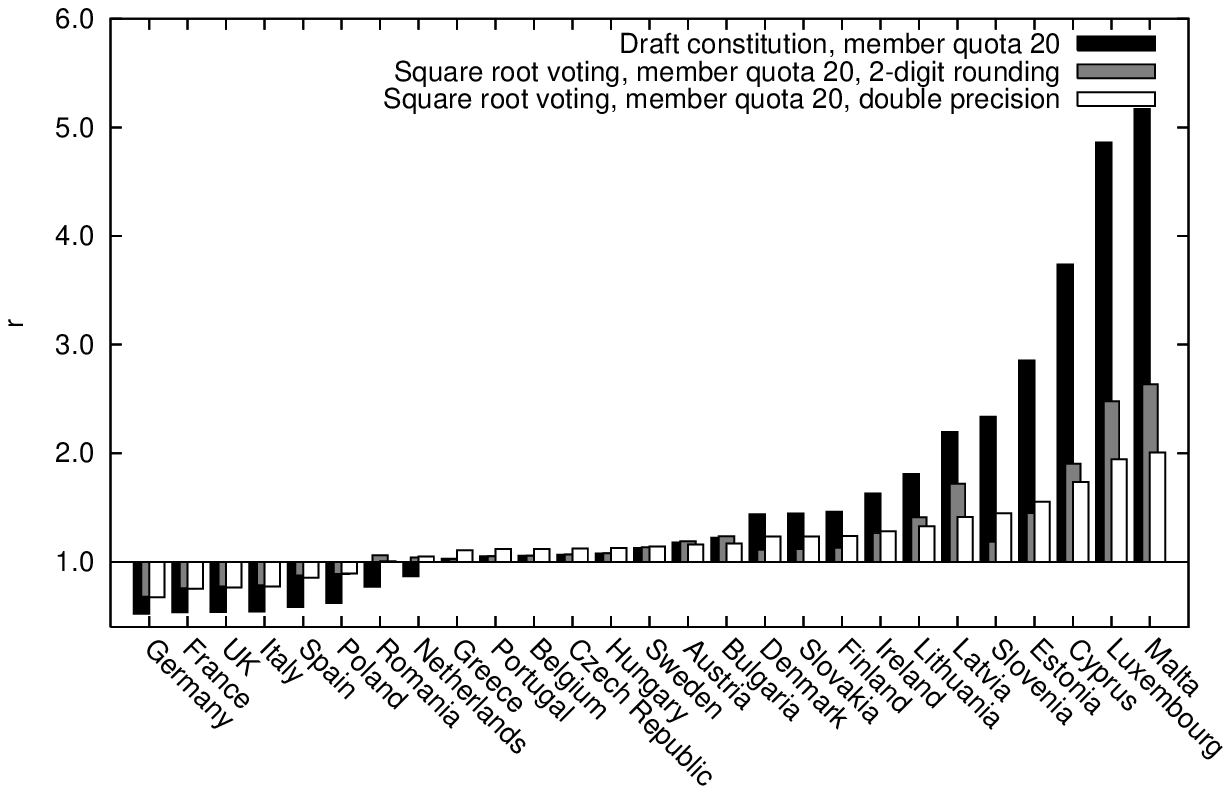}
\end{center}
\caption{The ratio $r$ for square root voting and the draft constitution with a member quota of 20.}
\label{f_quota20}
\end{figure}

For each individual country, the ratio $r$ defined in equation~\ref{e_ratiodef} can be calculated, and compared
with the corresponding ratio from the Nice Scheme and the draft constitution for the same member quota. The results
are shown in figures~\ref{f_quota14} to~\ref{f_quota20}. As expected, the member quota increases the voting power
of the smaller countries, but it can be seen that even with a member quota, the square
root voting system is fairer in the sense of Penrose.

A comparison between the square root voting system and the Nice/consti\-tution schemes should include the
efficiencies. Table~\ref{t_efficiency} shows the efficiency of square root voting with member quotas, and compares
them with the corresponding efficiencies in the other two schemes.

\begin{table}[ht]
\begin{center}
\begin{tabular}{|r|r|r|r|}
\hline
Member quota & \multicolumn{1}{|c|}{Square root} & \multicolumn{1}{c|}{Nice} & \multicolumn{1}{c|}{Draft constitution} \\
\hline
$14$ & $0.104526$ & $0.020256$ & --- \\
$15$ & $0.057132$ & --- & $0.128419$ \\
$18$ & $0.005479$ & $0.015809$ & --- \\
$20$ & $0.000904$ & --- & $0.007965$ \\
\hline
\end{tabular}
\end{center}
\caption{Efficiency of square root voting with member quota compared to Nice and the draft constitution.}
\label{t_efficiency}
\end{table}

It can be seen that for a member quota of 14, the square root voting system improves the efficiency compared to Nice
(the low efficiency of the Nice scheme was one of the main reasons for introducing a new voting scheme in the draft
constitution). For the other member quotas, the efficiency of the square root system is considerably lower than that
of the two existing schemes, which may be a reason not to adopt square root voting with these quotas in the form
described in this paper. The reason for the low efficiency is that the optimum range of quotas occurs at relatively
large values of $R$. However, the function $\sigma$ remains comparatively small even for values of $R$ that are somewhat
less than in the optimal range, and the low efficiency problem may be cured by choosing such a smaller value of $R$.
The price to pay is a somewhat increased value of $\sigma_{\rm min}$.

\section{Summary and conclusions}
\label{s_summary}

The concept of voting power according to Penrose, and its implementation by the Jagiellonian Compromise for the
Council of the European Union have been summarised.

It has been shown that the exact optimal quota for square root voting weights can be found analytically, but may
depend strongly on the rounding of the individual voting weights. It was demonstrated by two examples that rounding
may increase or decrease the deviation from fair voting in the sense of Penrose, and that rounding can be used
as a tool to find the optimal voting weights. The square root voting system has been analysed both on its own,
and with the inclusion of member quotas. For the 2007 population data from the member states of the EU, the
resulting voting weights, optimal quota ranges, and efficiencies are shown in tables~\ref{t_results1} and~\ref{t_results2}.

\begin{center}
\begin{longtable}{llr|r|r}
\cline{3-5}
 & & \multicolumn{3}{|c|}{Member quota} \\ 
\cline{3-5}
 & & \multicolumn{1}{|c}{None} & \multicolumn{1}{|c}{$14$} & \multicolumn{1}{|c|}{$15$} \\ 
\hline
\endfirsthead

\multicolumn{5}{l}{\tablename\ \thetable{} -- continued from previous page}\\
\cline{3-5}
 & & \multicolumn{3}{|c|}{Member quota} \\ 
\cline{3-5}
 & & \multicolumn{1}{|c}{None} & \multicolumn{1}{|c}{$14$} & \multicolumn{1}{|c|}{$15$}  \\ 
\hline
\endhead

\hline
\multicolumn{5}{r}{continued on next page}
\endfoot

\endlastfoot

\multicolumn{1}{|c}{ } & \multicolumn{1}{|l}{Germany} & \multicolumn{1}{|r}{$0.0944338$}& \multicolumn{1}{|r}{$0.094434$} & \multicolumn{1}{|r|}{$0.094434$}\\ 
\multicolumn{1}{|c}{ } & \multicolumn{1}{|l}{France} & \multicolumn{1}{|r}{$0.0828736$}& \multicolumn{1}{|r}{$0.082874$} & \multicolumn{1}{|r|}{$0.082874$} \\ 
\multicolumn{1}{|c}{ } & \multicolumn{1}{|l}{UK} & \multicolumn{1}{|r}{$0.0811605$}& \multicolumn{1}{|r}{$0.081160$} & \multicolumn{1}{|r|}{$0.081160$} \\ 
\multicolumn{1}{|c}{ } & \multicolumn{1}{|l}{Italy} & \multicolumn{1}{|r}{$0.0800400$}& \multicolumn{1}{|r}{$0.080040$} & \multicolumn{1}{|r|}{$0.080040$} \\ 
\multicolumn{1}{|c}{ } & \multicolumn{1}{|l}{Spain} & \multicolumn{1}{|r}{$0.0694152$}& \multicolumn{1}{|r}{$0.069415$} & \multicolumn{1}{|r|}{$0.069415$} \\ 
\multicolumn{1}{|c}{ } & \multicolumn{1}{|l}{Poland} & \multicolumn{1}{|r}{$0.0642697$}& \multicolumn{1}{|r}{$0.064270$} & \multicolumn{1}{|r|}{$0.064270$} \\ 
\multicolumn{1}{|c}{ } & \multicolumn{1}{|l}{Romania} & \multicolumn{1}{|r}{$0.0483364$}& \multicolumn{1}{|r}{$0.048336$} & \multicolumn{1}{|r|}{$0.048336$} \\ 
\multicolumn{1}{|c}{ } & \multicolumn{1}{|l}{Netherlands} & \multicolumn{1}{|r}{$0.0420982$}& \multicolumn{1}{|r}{$0.042098$} & \multicolumn{1}{|r|}{$0.042098$} \\ 
\multicolumn{1}{|c}{ } & \multicolumn{1}{|l}{Greece} & \multicolumn{1}{|r}{$0.0347891$}& \multicolumn{1}{|r}{$0.034789$} & \multicolumn{1}{|r|}{$0.034789$} \\ 
\multicolumn{1}{|c}{ } & \multicolumn{1}{|l}{Portugal} & \multicolumn{1}{|r}{$0.0338870$}& \multicolumn{1}{|r}{$0.033887$} & \multicolumn{1}{|r|}{$0.033887$} \\ 
\multicolumn{1}{|c}{ } & \multicolumn{1}{|l}{Belgium} & \multicolumn{1}{|r}{$0.0337465$}& \multicolumn{1}{|r}{$0.033746$} & \multicolumn{1}{|r|}{$0.033746$} \\ 
\multicolumn{1}{|c}{ } & \multicolumn{1}{|l}{Czech Republic} & \multicolumn{1}{|r}{$0.0333846$}& \multicolumn{1}{|r}{$0.033385$} & \multicolumn{1}{|r|}{$0.033385$} \\ 
\multicolumn{1}{|c}{ } & \multicolumn{1}{|l}{Hungary} & \multicolumn{1}{|r}{$0.0330205$}& \multicolumn{1}{|r}{$0.033021$} & \multicolumn{1}{|r|}{$0.033021$} \\ 
\multicolumn{1}{|c}{$v_i$} & \multicolumn{1}{|l}{Sweden} & \multicolumn{1}{|r}{$0.0314221$}& \multicolumn{1}{|r}{$0.031422$} & \multicolumn{1}{|r|}{$0.031422$} \\ 
\multicolumn{1}{|c}{ } & \multicolumn{1}{|l}{Austria} & \multicolumn{1}{|r}{$0.0299854$}& \multicolumn{1}{|r}{$0.029985$} & \multicolumn{1}{|r|}{$0.029985$} \\ 
\multicolumn{1}{|c}{ } & \multicolumn{1}{|l}{Bulgaria} & \multicolumn{1}{|r}{$0.0288442$}& \multicolumn{1}{|r}{$0.028844$} & \multicolumn{1}{|r|}{$0.028844$} \\ 
\multicolumn{1}{|c}{ } & \multicolumn{1}{|l}{Denmark} & \multicolumn{1}{|r}{$0.0242930$}& \multicolumn{1}{|r}{$0.024293$} & \multicolumn{1}{|r|}{$0.024293$} \\ 
\multicolumn{1}{|c}{ } & \multicolumn{1}{|l}{Slovakia} & \multicolumn{1}{|r}{$0.0241735$}& \multicolumn{1}{|r}{$0.024173$} & \multicolumn{1}{|r|}{$0.024173$} \\ 
\multicolumn{1}{|c}{ } & \multicolumn{1}{|l}{Finland} & \multicolumn{1}{|r}{$0.0239106$}& \multicolumn{1}{|r}{$0.023911$} & \multicolumn{1}{|r|}{$0.023911$} \\ 
\multicolumn{1}{|c}{ } & \multicolumn{1}{|l}{Ireland} & \multicolumn{1}{|r}{$0.0213545$}& \multicolumn{1}{|r}{$0.021354$} & \multicolumn{1}{|r|}{$0.021354$} \\ 
\multicolumn{1}{|c}{ } & \multicolumn{1}{|l}{Lithuania} & \multicolumn{1}{|r}{$0.0191501$}& \multicolumn{1}{|r}{$0.019150$} & \multicolumn{1}{|r|}{$0.019150$} \\ 
\multicolumn{1}{|c}{ } & \multicolumn{1}{|l}{Latvia} & \multicolumn{1}{|r}{$0.0157214$}& \multicolumn{1}{|r}{$0.015721$} & \multicolumn{1}{|r|}{$0.015721$} \\ 
\multicolumn{1}{|c}{ } & \multicolumn{1}{|l}{Slovenia} & \multicolumn{1}{|r}{$0.0147583$}& \multicolumn{1}{|r}{$0.014758$} & \multicolumn{1}{|r|}{$0.014758$} \\ 
\multicolumn{1}{|c}{ } & \multicolumn{1}{|l}{Estonia} & \multicolumn{1}{|r}{$0.0120598$}& \multicolumn{1}{|r}{$0.012060$} & \multicolumn{1}{|r|}{$0.012060$} \\ 
\multicolumn{1}{|c}{ } & \multicolumn{1}{|l}{Cyprus} & \multicolumn{1}{|r}{$0.0091841$}& \multicolumn{1}{|r}{$0.009184$} & \multicolumn{1}{|r|}{$0.009184$} \\ 
\multicolumn{1}{|c}{ } & \multicolumn{1}{|l}{Luxembourg} & \multicolumn{1}{|r}{$0.0070557$}& \multicolumn{1}{|r}{$0.007056$} & \multicolumn{1}{|r|}{$0.007056$} \\ 
\multicolumn{1}{|c}{ } & \multicolumn{1}{|l}{Malta} & \multicolumn{1}{|r}{$0.0066324$}& \multicolumn{1}{|r}{$0.006632$} & \multicolumn{1}{|r|}{$0.006632$} \\ 
\hline
\multicolumn{2}{|l}{Quota $>$} & \multicolumn{1}{|r|}{$0.6149670$}& $0.6146660$ & \multicolumn{1}{r|}{$0.682884$} \\ 
\multicolumn{2}{|l}{Quota $\leq$} & \multicolumn{1}{|r|}{$0.6149671$} & $0.6146661$ & \multicolumn{1}{r|}{$0.682885$} \\ 
\hline
\multicolumn{2}{|l}{Efficiency} & \multicolumn{1}{|r|}{$0.1644$} & $0.1045$ & \multicolumn{1}{r|}{$0.05713$} \\ 
\hline
\caption{Summary of optimal voting weights, optimal quota ranges, and efficiencies in square root voting with different member quotas.}
\label{t_results1}
\end{longtable}

\end{center}

\begin{center}
\begin{longtable}{llr|r}
\cline{3-4}
 & & \multicolumn{2}{|c|}{Member quota} \\ 
\cline{3-4}
 & & \multicolumn{1}{|c}{$18$} & \multicolumn{1}{|c|}{$20$} \\ 
\hline
\endfirsthead

\multicolumn{4}{l}{\tablename\ \thetable{} -- continued from previous page}\\
\cline{3-4}
 & & \multicolumn{2}{|c|}{Member quota} \\ 
\cline{3-4}
 & & \multicolumn{1}{|c}{$18$} & \multicolumn{1}{|c|}{$20$} \\ 
\hline
\endhead

\hline
\multicolumn{4}{r}{continued on next page}
\endfoot

\endlastfoot

\multicolumn{1}{|c}{ } & \multicolumn{1}{|l}{Germany} & \multicolumn{1}{|r}{$0.094434$} & \multicolumn{1}{|r|}{$0.09$} \\ 
\multicolumn{1}{|c}{ } & \multicolumn{1}{|l}{France} & \multicolumn{1}{|r}{$0.082874$} & \multicolumn{1}{|r|}{$0.08$} \\ 
\multicolumn{1}{|c}{ } & \multicolumn{1}{|l}{UK} & \multicolumn{1}{|r}{$0.081160$} & \multicolumn{1}{|r|}{$0.08$} \\ 
\multicolumn{1}{|c}{ } & \multicolumn{1}{|l}{Italy} & \multicolumn{1}{|r}{$0.080040$} & \multicolumn{1}{|r|}{$0.08$} \\ 
\multicolumn{1}{|c}{ } & \multicolumn{1}{|l}{Spain} & \multicolumn{1}{|r}{$0.069415$} & \multicolumn{1}{|r|}{$0.07$} \\ 
\multicolumn{1}{|c}{ } & \multicolumn{1}{|l}{Poland} & \multicolumn{1}{|r}{$0.064270$} & \multicolumn{1}{|r|}{$0.06$} \\ 
\multicolumn{1}{|c}{ } & \multicolumn{1}{|l}{Romania} & \multicolumn{1}{|r}{$0.048336$} & \multicolumn{1}{|r|}{$0.05$} \\ 
\multicolumn{1}{|c}{ } & \multicolumn{1}{|l}{Netherlands} & \multicolumn{1}{|r}{$0.042098$} & \multicolumn{1}{|r|}{$0.04$} \\ 
\multicolumn{1}{|c}{ } & \multicolumn{1}{|l}{Greece} & \multicolumn{1}{|r}{$0.034789$} & \multicolumn{1}{|r|}{$0.03$} \\ 
\multicolumn{1}{|c}{ } & \multicolumn{1}{|l}{Portugal} & \multicolumn{1}{|r}{$0.033887$} & \multicolumn{1}{|r|}{$0.03$} \\ 
\multicolumn{1}{|c}{ } & \multicolumn{1}{|l}{Belgium} & \multicolumn{1}{|r}{$0.033746$} & \multicolumn{1}{|r|}{$0.03$} \\ 
\multicolumn{1}{|c}{ } & \multicolumn{1}{|l}{Czech Republic} & \multicolumn{1}{|r}{$0.033385$} & \multicolumn{1}{|r|}{$0.03$} \\ 
\multicolumn{1}{|c}{ } & \multicolumn{1}{|l}{Hungary} & \multicolumn{1}{|r}{$0.033021$} & \multicolumn{1}{|r|}{$0.03$} \\ 
\multicolumn{1}{|c}{$v_i$} & \multicolumn{1}{|l}{Sweden} & \multicolumn{1}{|r}{$0.031422$} & \multicolumn{1}{|r|}{$0.03$} \\ 
\multicolumn{1}{|c}{ } & \multicolumn{1}{|l}{Austria} & \multicolumn{1}{|r}{$0.029985$} & \multicolumn{1}{|r|}{$0.03$} \\ 
\multicolumn{1}{|c}{ } & \multicolumn{1}{|l}{Bulgaria} & \multicolumn{1}{|r}{$0.028844$} & \multicolumn{1}{|r|}{$0.03$} \\ 
\multicolumn{1}{|c}{ } & \multicolumn{1}{|l}{Denmark} & \multicolumn{1}{|r}{$0.024293$} & \multicolumn{1}{|r|}{$0.02$} \\ 
\multicolumn{1}{|c}{ } & \multicolumn{1}{|l}{Slovakia} & \multicolumn{1}{|r}{$0.024173$} & \multicolumn{1}{|r|}{$0.02$} \\ 
\multicolumn{1}{|c}{ } & \multicolumn{1}{|l}{Finland} & \multicolumn{1}{|r}{$0.023911$} & \multicolumn{1}{|r|}{$0.02$} \\ 
\multicolumn{1}{|c}{ } & \multicolumn{1}{|l}{Ireland} & \multicolumn{1}{|r}{$0.021354$} & \multicolumn{1}{|r|}{$0.02$} \\ 
\multicolumn{1}{|c}{ } & \multicolumn{1}{|l}{Lithuania} & \multicolumn{1}{|r}{$0.019150$} & \multicolumn{1}{|r|}{$0.02$} \\ 
\multicolumn{1}{|c}{ } & \multicolumn{1}{|l}{Latvia} & \multicolumn{1}{|r}{$0.015721$} & \multicolumn{1}{|r|}{$0.02$} \\ 
\multicolumn{1}{|c}{ } & \multicolumn{1}{|l}{Slovenia} & \multicolumn{1}{|r}{$0.014758$} & \multicolumn{1}{|r|}{$0.01$} \\ 
\multicolumn{1}{|c}{ } & \multicolumn{1}{|l}{Estonia} & \multicolumn{1}{|r}{$0.012060$} & \multicolumn{1}{|r|}{$0.01$} \\ 
\multicolumn{1}{|c}{ } & \multicolumn{1}{|l}{Cyprus} & \multicolumn{1}{|r}{$0.009184$} & \multicolumn{1}{|r|}{$0.01$} \\ 
\multicolumn{1}{|c}{ } & \multicolumn{1}{|l}{Luxembourg} & \multicolumn{1}{|r}{$0.007056$} & \multicolumn{1}{|r|}{$0.01$} \\ 
\multicolumn{1}{|c}{ } & \multicolumn{1}{|l}{Malta} & \multicolumn{1}{|r}{$0.006632$} & \multicolumn{1}{|r|}{$0.01$} \\ 
\hline
\multicolumn{2}{|l}{Quota $>$} & \multicolumn{1}{|r|}{$0.784222$} & \multicolumn{1}{r|}{$0.80$} \\ 
\multicolumn{2}{|l}{Quota $\leq$} & \multicolumn{1}{|r|}{$0.784223$} & \multicolumn{1}{r|}{$0.81$} \\ 
\hline
\multicolumn{2}{|l}{Efficiency} & \multicolumn{1}{|r|}{$0.005479$} & \multicolumn{1}{r|}{$0.000904$} \\ 
\hline
\caption{Summary of optimal voting weights, optimal quota ranges, and efficiencies in square root voting with different member quotas.}
\label{t_results2}
\end{longtable}

\end{center}

It has been shown that these voting schemes have a smaller deviation from an ideal square root voting power system
than both the Nice voting scheme and the draft constitution voting scheme. However, for larger member quotas, the
square root voting schemes minimise this deviation for rather large vote quotas, resulting in very small efficiency.
To use square root voting with large member quotas, a compromise between the desire for a ``fair'' voting scheme and
an efficient system would have to be found.

\bibliography{paper}

\begin{thebibliography}{33}
\expandafter\ifx\csname natexlab\endcsname\relax\def\natexlab#1{#1}\fi
\expandafter\ifx\csname url\endcsname\relax
  \def\url#1{\texttt{#1}}\fi
\expandafter\ifx\csname urlprefix\endcsname\relax\def\urlprefix{URL }\fi

\bibitem[{Albert(2003)}]{Albert1}
Albert, M., 2003. The voting power approach: Measurement without theory.
  European Union Politics 4, 351--366.

\bibitem[{Bafumi et~al.(2004)Bafumi, Gelman, and Katz}]{bafumi1}
Bafumi, J., Gelman, A., Katz, J.~N., 2004. Standard voting power indexes do not
  work: An empirical analysis. B.J.Pol.S. 34, 657--674.

\bibitem[{Baldwin and Widgren(2003{\natexlab{a}})}]{Baldwin1}
Baldwin, R., Widgren, M., 2003{\natexlab{a}}. Decision making and the
  {C}onstitutional {T}reaty: {W}ill the {IGC} discard {G}iscard?, {C}entre for
  {E}uropean {P}olicy {S}tudies, {P}olicy {B}riefs {N}o 37.

\bibitem[{Baldwin and Widgren(2003{\natexlab{b}})}]{Baldwin2}
Baldwin, R., Widgren, M., 2003{\natexlab{b}}. A study of the {C}onstitutional
  {T}reaty's voting reform dilemma. {C}entre for {E}uropean {P}olicy {S}tudies,
  {P}olicy {B}riefs {N}o 37.

\bibitem[{Baldwin and Widgren(2004)}]{Baldwin3}
Baldwin, R., Widgren, M., 2004. Winners and losers under various dual majority
  rules for the {EU}'s {C}ouncil of {M}inisters, {CEPS} {P}olicy {B}rief {N}o
  50.

\bibitem[{Banzhaf(1965)}]{banzhaf1}
Banzhaf, J.~F., 1965. Weighted voting does not work: A mathematical analysis.
  Rutgers Law Review 19, 317--343.

\bibitem[{Banzhaf(1968)}]{banzhaf2}
Banzhaf, J.~F., 1968. One man, 3312 votes: {A} mathematical analysis of the
  {E}lectoral {C}ollege. Villanova Law Review 13, 303--332.

\bibitem[{Barber{\`a} and Jackson(2004)}]{Barbera1}
Barber{\`a}, S., Jackson, M.~O., 2004. On the weights of nations: {A}ssigning
  voting weights in a heterogeneous union.
\newline\urlprefix\url{http://citeseer.ist.psu.edu/605197.html}

\bibitem[{Coleman(1971)}]{Coleman1}
Coleman, J.~S., 1971. Control of collectivities and the power of a collectivity
  to act, in {B.~Lieberman} (ed), {Social Choice}, {New York}.

\bibitem[{{EUobserver}(2007)}]{euobserver}
{EUobserver}, 2007. E{U}observer.com, {P}oland to fight for square root law in
  {EU} {T}reaty, 29 {M}arch 2007.
\newline\urlprefix\url{http://www.euobserver.com}

\bibitem[{{European Statistics Office}(2007)}]{EuroStat}
{European Statistics Office}, 2007. eurostat.
\newline\urlprefix\url{http://ec.europa.eu/eurostat}

\bibitem[{{European Union}(2001)}]{NiceTreaty}
{European Union}, 2001. Treaty of {N}ice, {O}fficial {J}ournal of the
  {E}uropean {U}nion {C} 80 (2001).

\bibitem[{{European Union}(2003)}]{2003AccTreaty}
{European Union}, 2003. Treaty between the {M}ember {S}tates of the {E}uropean
  {U}nion and the {C}zech {R}epublic, the {R}epublic of {E}stonia, the
  {R}epublic of {C}yprus, the {R}epublic of {L}atvia, the {R}epublic of
  {L}ithuania, the {R}epublic of {H}ungary, the {R}epublic of {M}alta, the
  {R}epublic of {P}oland, the {R}epublic of {S}lovenia, the {S}lovak
  {R}epublic, concerning the accession of the {C}zech {R}epublic, the
  {R}epublic of {E}stonia, the {R}epublic of {C}yprus, the {R}epublic of
  {L}atvia, the {R}epublic of {L}ithuania, the {R}epublic of {H}ungary, the
  {R}epublic of {M}alta, the {R}epublic of {P}oland, the {R}epublic of
  {S}lovenia and the {S}lovak {R}epublic to the {E}uropean {U}nion, {O}fficial
  {J}ournal of the {E}uropean {U}nion {L} 236 (2003).

\bibitem[{{European Union}(2004)}]{DraftConstitution}
{European Union}, 2004. Treaty establishing a {C}onstitution for {E}urope,
  {O}fficial {J}ournal of the {E}uropean {U}nion {C} 310 (2004).

\bibitem[{{European Union}(2005)}]{2005AccTreaty}
{European Union}, 2005. Treaty between the {M}ember {S}tates of the {E}uropean
  {U}nion and the {R}epublic of {B}ulgaria and {R}omania, concerning the
  accession of the {R}epublic of {B}ulgaria and {R}omania to the {E}uropean
  {U}nion, {O}fficial {J}ournal of the {E}uropean {U}nion {L} 157 (2005).

\bibitem[{{European Union}(2007)}]{ReformTreaty}
{European Union}, 2007. Conference of the {R}epresentatives of the
  {G}overnments of the {M}ember {S}tates, {D}raft {T}reaty {A}mending the
  {T}reaty on {E}uropean {U}nion and the {T}reaty {E}stablishing the {E}uropean
  {C}ommunity, {CIG} 1/1/07 {REV} 1.
\newline\urlprefix\url{www.consilium.europa.eu}

\bibitem[{Felderer et~al.(2003)Felderer, Paterson, and
  Sil{\'a}rsky}]{Felderer1}
Felderer, B., Paterson, I., Sil{\'a}rsky, P., 2003. The double majority implies
  a massive transfer of power to the large member states -- is this intended?
  {S}hort {P}olicy {P}aper for {EU} {C}onvention {F}orum, {I}nstitute for
  {A}dvanced {S}tudies, {V}ienna.

\bibitem[{Felsenthal and Machover(2004)}]{Felsenthal1}
Felsenthal, D.~S., Machover, M., 2004. Analysis of {QM} rules in the draft
  constitution for {E}urope proposed by the {European Convention}, 2003. Social
  Choice and Welfare 23~(1), 1--20.

\bibitem[{Garret and Tsebelis(2001)}]{Garrett1}
Garret, G., Tsebelis, G., 2001. Even more reasons to resist the temptation of
  power indices in the {EU}. Journal of Theoretical Politics 13~(1), 99--105.

\bibitem[{Hoare(1961)}]{hoare1}
Hoare, C.~A.~R., 1961. Quicksort: {A}lgorithm 64. Comm.~ACM 4, 321--322.

\bibitem[{Hoare(1962)}]{hoare2}
Hoare, C.~A.~R., 1962. Quicksort. Computer J. 5, 10--15.

\bibitem[{Hosli and Machover(2004)}]{Hosli1}
Hosli, M.~O., Machover, M., 2004. The {N}ice {T}reaty and voting rules in the
  {C}ouncil: {A} reply to {M}oberg (2002). Journal of Common Market Studies
  42~(3), 497--521.

\bibitem[{Leech(2002)}]{Leech1}
Leech, D., 2002. Designing the voting system for the {EU} {C}ouncil of
  {M}inisters. Public Choice 113, 437--464.

\bibitem[{Leech(2003)}]{Leech2}
Leech, D., 2003. The utility of the voting power approach, {CSGR} {W}orking
  {P}aper {N}o.~118/03.
\newline\urlprefix\url{http://www2.warwick.ac.uk/fac/soc/\\
  csgr/research/workingpapers/2003/wp11803.pdf}

\bibitem[{Nurmi and Meskanen(1999)}]{Nurmi1}
Nurmi, H., Meskanen, R., 1999. A priori power measures and the institutions of
  the {European Union}. European Journal of Political Research 35, 161--179.

\bibitem[{Penrose(1946)}]{penrose1}
Penrose, L.~S., 1946. The elementary statistics of majority voting. Journal of
  the Royal Statistical Society 109~(1), 53--57.

\bibitem[{Plechanovov{\'a}(2003)}]{Plechanovova1}
Plechanovov{\'a}, B., 2003. The {T}reaty of {N}ice and the distribution of
  votes in the {C}ouncil -- {V}oting power consequences for the {EU} after the
  oncoming enlargement. European Integration online Papers 7~(6).
\newline\urlprefix\url{http://eiop.or.at/eiop/pdf/2003-006.pdf}

\bibitem[{Press et~al.(1993)Press, Flannery, Teukolski, and
  Vetterling}]{NumericalRecipes}
Press, W.~H., Flannery, B.~P., Teukolski, S.~A., Vetterling, W.~T., 1993.
  Numerical Recipes in C, 2nd Edition. Cambridge University Press.

\bibitem[{{Scientists for a Democratic Europe}(2007)}]{OpenLetter}
{Scientists for a Democratic Europe}, 2007. Letter to the governments of the
  {EU} member states.
\newline\urlprefix\url{http://www.ruhr-uni-bochum.de/\\
  mathphys/politik/eu/open-letter.htm}

\bibitem[{Shapley and Shubik(1954)}]{Shapley1}
Shapley, L.~S., Shubik, M., 1954. A method for evaluating the distribution of
  power in a committee system. American Political Science Review 48, 787--792.

\bibitem[{Zyczkowski and Slomczynski(2004)}]{slomczynski1}
Zyczkowski, K., Slomczynski, W., 2004. Voting in the {E}uropean {U}nion: The
  square root system of {P}enrose and a critical point.
\newline\urlprefix\url{http://www.citebase.org/\\abstract?id=oai:arXiv.org:con%
d-mat/0405396}

\bibitem[{Zyczkowski and Slomczynski(2006)}]{slomczynski3}
Zyczkowski, K., Slomczynski, W., 2006. Penrose voting system and optimal quota.
  Physica Acta Polonica B 37~(11), 3133.

\bibitem[{Zyczkowski et~al.(2006)Zyczkowski, Slomczynski, and
  Zastawniak}]{slomczynski2}
Zyczkowski, K., Slomczynski, W., Zastawniak, T., 2006. Physics for fairer
  voting. Physics World 19~(3), 35--37.

\end{thebibliography}

\end{document}